\documentclass[10pt,leqno]{article}

\usepackage{lmodern,color}
\usepackage[french]{babel}
\usepackage{amsmath,amssymb,amscd}
\usepackage{a4wide}
\usepackage[utf8]{inputenc}
\usepackage[all]{xy} 
\usepackage[T1]{fontenc}
\usepackage{hyperref}
\usepackage{chngcntr}
\usepackage{comment}

\title{Un caractère relatif pondéré}
\author{Pierre-Henri Chaudouard}
\date{ version du \today}

\newenvironment{paragr}[1][]{\refstepcounter{subsubsection} \noindent \textbf{\thesubsubsection . \ #1}}{\medskip}

\newenvironment{theoreme}{ \medskip\refstepcounter{theo}  \noindent\textbf{Th\'eor\`eme \thetheo}. ---\em}{\em \medskip}
\newenvironment{proposition}{\medskip\refstepcounter{theo}   \noindent\textbf{Proposition \thetheo}. ---\em}{\em\medskip}
\newenvironment{corollaire}{\medskip\refstepcounter{theo}  \noindent\textbf{Corollaire \thetheo}. ---\em}{\em\medskip}

\newenvironment{lemme}{\medskip\refstepcounter{theo}   \noindent\textbf{Lemme \thetheo}. ---\em}{\em\medskip}

\newenvironment{preuve}[1][]{\noindent \textbf{Démonstration.} #1 --- }{\hfill
  \ensuremath{\square} \medskip}

\newenvironment{remarque}{\medskip\refstepcounter{theo}  \noindent\textbf{Remarque \thetheo}. ---}{\medskip}

\DeclareMathOperator{\vol}{vol}

\DeclareMathOperator{\Hom}{Hom}

\DeclareMathOperator{\Res}{Res}

\DeclareMathOperator{\Ker}{Ker}

\DeclareMathOperator{\vect}{vect}

\DeclareMathOperator{\Sym}{Sym}

\newcommand{\ZZ}{\mathbb{Z}}

\newcommand{\RR}{\mathbb{R}}
\newcommand{\AAA}{\mathbb{A}}
\newcommand{\CC}{\mathbb{C}}

\newcommand{\ga}{\gamma}

\newcommand{\oc}{\mathcal{O}}
\newcommand{\rc}{\mathcal{R}}
\newcommand{\Sc}{\mathcal{S}}

\newcommand{\ec}{\mathcal{E}}

\newcommand{\mc}{\mathcal{M}}

\newcommand{\fc}{\mathcal{F}}

\newcommand{\Ac}{\mathcal{A}}
\newcommand{\bc}{\mathcal{B}}

\newcommand{\Jc}{\mathcal{J}}

\newcommand{\Gc}{\widehat{G}}

\newcommand{\ago}{\mathfrak{a}}

\newcommand{\Xgo}{\mathfrak{X}}

\newcommand{\al}{\alpha}
\newcommand{\be}{\beta}

\newcommand{\La}{\Lambda}
\newcommand{\la}{\lambda}

\newcommand{\back}{\backslash}

\newcommand{\bg}{\langle}
\newcommand{\bd}{\rangle}

\newcommand{\eps}{\varepsilon}

\renewcommand{\leq}{\leqslant}
\renewcommand{\geq}{\geqslant}

\begin{document}
\counterwithin{equation}{subsubsection}

\maketitle

\begin{abstract}
  Soit $p\geq 1$.   L'espace symétrique $S=GL(2p+1)/GL(p+1)\times GL(p)$, sur un corps de nombres, n'est pas cuspidal au sens où son spectre automorphe ne comporte aucune représentation cuspidale de $GL(2p+1)$. Dans cet article, on calcule  la décomposition spectrale  de sa partie relativement cuspidale: par définition, c'est  l'induite de la partie cuspidale de l'espace symétrique $(GL(1)\times GL(2p)) /  (GL(1)\times GL(p)\times GL(p))$. Comme application, on obtient l'expression de la contribution de cette partie relativement cuspidale à la formule des traces de Guo-Jacquet (établie par H. Li et l'auteur) sous la forme d'un caractère relatif pondéré.
\end{abstract}

{\selectlanguage{english}
\begin{abstract}
  Let $p\geq 1$. The symmetric space $S=GL(2p+1)/GL(p+1)\times GL(p)$ (over a number field) is not cuspidal in the sense that its automorphic spectrum does not contain any cuspidal representation of $GL(2p+1)$. In this article, we compute the spectral decomposition of its relatively cuspidal part: this  is, by definition, the part of the spectrum that is   induced from the cuspidal part of the symmetric space $(GL(1)\times GL(2p)) / (GL(1)\times GL(p)\times GL(p))$. As an application, we obtain the expression of the contribution of this relatively cuspidal part to the Guo-Jacquet trace formula (established by H. Li and the author)  in terms of a weighted relative character.
\end{abstract} }

\tableofcontents

\section{Introduction}

\subsection{Le spectre automorphe des espaces symétriques}

\begin{paragr}
  Commençons par un problème soulevé, entre autres, par Jacquet et ses collaborateurs, cf.  \cite{FriJa,Jac-Edin}. Soit $F$ un corps de nombres et $\AAA$ l'anneau des adèles de $F$. Soit $G$ un groupe algébrique affine, réductif, connexe, défini sur $F$ et $\theta$ un automorphisme involutif du $F$-groupe algébrique $G$. Soit
  \begin{align*}
    S_0=\{g \in G\mid \theta(g)=g^{-1}\}.
  \end{align*}
  Le groupe $G$ agit sur $S_0$ par $\theta$-conjugaison qui est donnée, pour tous $g\in G$ et $s\in S_0$, par  $g\cdot s=g s\theta(g)^{-1}\in S_0$. Soit $\xi_0\in S_0(F)$ et $S$ l'espace symétrique défini par la $G$-orbite de $S_0$. Soit $\Phi$ une fonction lisse à support compact sur $S(\AAA)$ (plus généralement une fonction de Schwartz). On peut alors former la série
  \begin{align*}
    K_\Phi(g)=\sum_{\xi \in S(F)} \Phi(g^{-1}\cdot \xi). 
  \end{align*}
  On obtient ainsi une fonction sur $G(\AAA)$ invariante à gauche par $G(F)$ et $Z_G^\theta(\AAA)$, où $Z_G^\theta$ est le sous-groupe du centre de $G$ formé des éléments fixes sous $\theta$. Le groupe $G(\AAA)$ agit par translation à droite sur l'espace engendré par les fonctions $K_\Phi$ lorsque $\Phi$ décrit l'espace de fonctions test considérées. Un problème fondamental est alors d'obtenir la décomposition  spectrale de cet espace. Signalons  au moins deux cas où ce problème a été résolu: il s'agit du cas où $G=G_1\times G_1$ et $\theta$ est l'involution qui échange les deux facteurs. Dans ce cas, $S_0$ s'identifie à $G_1$ muni de l'action (transitive) de $G$ donnée par $(x,y)\cdot g=x g y^{-1}$. La fonction $K_\Phi$ est alors le \og noyau automorphe \fg{}, noyau de l'action par convolution de $\Phi$  sur  $L^2(G_1(F)\back G_1(\AAA))$, dont le développement spectral repose sur la décomposition de Langlands de $L^2(G_1(F)\back G_1(\AAA))$, cf. par exemple \cite[éq. (12.7) et (12.8)]{Ar-cours}. Le second cas est celui  considéré dans \cite{chaudouardsymmetric}: il s'agit du groupe $G=\Res_{E/F}(GL(n,E))$ muni de l'involution galoisienne, pour un entier $n\geq 1$ et $E$ une extension quadratique de $F$.
\end{paragr}

\begin{paragr} Revenons au  cas d'un couple $(G,S)$ général. Un problème plus immédiatement abordable est la détermination de la composante cuspidale de $K_\Phi$. Pour simplifier la discussion, on suppose que le groupe $G(F)$ agit transitivement sur $S(F)$. Soit $H$ le stabilisateur de $\xi_0$. Soit $(\pi,V_\pi)$ une représentation automorphe cuspidale de $G(\AAA)$ de caractère central unitaire et trivial sur $Z_G^\theta(\AAA)$. Soit $\phi$ une forme automorphe dans l'espace $V_\pi$ de $\pi$. Selon le calcul de Friedberg-Jacquet, \cite[section 1]{FriJa}, on a
  \begin{align*}
    \int_{Z_G^\theta(\AAA)G(F)\back G(\AAA)    } K_\Phi(g) \overline{\phi(g)}\, dg=\int_{  H(\AAA)  \back G(\AAA)   } \Phi(g^{-1}\cdot \xi_0) \overline{\int_{Z_G^\theta(\AAA)H(F)\back H(\AAA)  }  \phi(hg)\, dh} \, dg.
  \end{align*}
  Il s'ensuit que $K_\Phi$ a une composante non triviale sur $V_\pi$ si et seulement si $V_\pi$ est $H$-distingué au sens où la $H$-période, définie comme la forme linéaire
  \begin{align}\label{eq:intro-H per}
    \phi\in V_\pi \mapsto \int_{Z_G^\theta(\AAA)H(F)\back H(\AAA)  }  \phi(h)\, dh,
  \end{align}
  n'est pas identiquement nulle. Notons que l'intégrale est bien définie: le groupe $H$ est réductif,  le quotient $Z_G^\theta(\AAA)H(F)\back H(\AAA)$ s'identifie à un fermé de $Z_G(\AAA)G(F)\back G(\AAA)$  et la fonction cuspidale $\phi$ est à décroissance rapide sur ce dernier. 
\end{paragr}

\begin{paragr} \label{S:intro-gp dual} On s'attend à ce que les représentations automorphes qui apparaissent dans la décomposition spectrale de $K_\Phi$ aient des propriétés remarquables vis-à-vis de la fonctorialité de Langlands. Lorsque $G$ est déployé, Sakellaridis et Venkatesh ont défini dans $\cite{SaVen}$, cf. aussi \cite{Sak-spherical} et \cite{takeda2023}, un groupe dual $\Gc_S$ muni d'un morphisme $\rho$ de $\Gc_S\times SL(2,\CC)$ vers le groupe dual $\Gc$ de $G$ de sorte que les paramètres de Langlands, associés conjecturalement aux représentations automorphes de $G$ qui apparaissent dans le spectre des fonctions $K_\Phi$, devraient se factoriser par $\rho$. 
\end{paragr}

\begin{paragr} \label{S:intro-GJ} Pour tout $n\geq 1$, soit $I_n$ la matrice identité de taille $n$. Désormais, nous  fixons $n\geq 1$ et nous  nous focalisons sur  le cas du groupe $G=GL(n)$ sur le corps $F$ et de l'involution donnée par l'automorphisme intérieur associé à $\theta=
  \begin{pmatrix}
    I_p & 0 \\ 0 & -I_q
  \end{pmatrix}$ pour des entiers $n\geq 2$ et $p,q\geq 0$ tels que  $n=p+q$. Dans ce cas,
  \begin{align*}
    S_0=\{ g\in G \mid (g\theta)^2=I_n\}.
  \end{align*}
  Soit $S\subset S_0$ la $G$-orbite (pour la $\theta$-conjugaison) de $I_n$.  Le stabilisateur de $I_n$ pour cette action est noté $H$ et est naturellement isomorphe à $GL(p)\times GL(q)$. Notons que $S(F)$ est aussi la $G(F)$-orbite de $I_n$. La fonction test $\Phi$ sur $S(\AAA)$ se déduit alors,  par intégration le long des fibres de l'application $g\in G(\AAA)\mapsto g\theta g^{-1}\theta\in S(\AAA)$,  d'une fonction  $f$ dans la classe de Schwartz $\Sc(G(\AAA))$ de $G(\AAA)$. On a donc la relation, pour tout $g\in G(\AAA)$,
  \begin{align*}
    \Phi(g\theta g^{-1}\theta)=\int_{H(\AAA)  }  f( g h)\, dh
  \end{align*}
  pour une mesure de Haar $dh$ fixée sur $H(\AAA)$. En introduisant le \og noyau automorphe \fg{}
  \begin{align*}
    K_f(g_1,g_2)=\sum_{\ga \in G(F)} f(g_1^{-1}\ga g_2)
  \end{align*}
  pour tous $g_1,g_2\in G(\AAA)$, on voit qu'on  a alors, pour tout $g\in G(\AAA)$,
  \begin{align}\label{eq:intro-KPhi}
    K_\Phi(g)= \int_{[H]}K_f(g, h)\, dh
  \end{align}
où $[H]=H(F)\back H(\AAA)$ et $dh$ est maintenant la mesure quotient de la mesure de Haar ci-dessus par la mesure de comptage sur $H(F)$.
  Supposons, pour fixer les idées, $p\geq q$. Dans ce cas, le groupe dual $\Gc_S$ est le groupe $Sp(2q,\CC)$ (groupe des automorphismes d'une forme symplectique non dégénérée sur $\CC^{2q}$). Si $p=q$ (et donc $n=2q$), le morphisme
  \begin{align*}
      \rho:\Gc_S\times SL(2,\CC)\to \Gc=GL(n)
  \end{align*}
   est la représentation naturelle sur le premier facteur et est trivial sur le second. Supposons $p>q$.  Soit $\Sym^{p-q-1}$ la représentation de $SL(2,\CC)$ de dimension $p-q$ donnée par la puissance symétrique de degré $p-q-1$ de la représentation standard. En sommant la représentation naturelle de $Sp(2q,\CC)$ et la représentation  $\Sym^{p-q-1}$, on obtient un morphisme   $Sp(2q,\CC)\times SL(2,\CC)\to GL(2q,\CC)\times GL(p-q,\CC)$. En identifiant $GL(2q,\CC)\times GL(p-q,\CC)$ à un sous-groupe de Levi de $GL(2n,\CC)$, on obtient alors le morphisme
   \begin{align*}
        \rho:\Gc_S\times SL(2,\CC)\to  \Gc
   \end{align*}
pour $p>q$. Dans ce cas, les conjectures du § \ref{S:intro-gp dual} prédisent alors que la composante cuspidale de $K_\Phi$ est nulle ou, ce qui revient au même, que les représentations cuspidales de $GL(n)$ ne sont jamais $H$-distinguées. C'est exactement ce que prouvent Jacquet et Friedberg, cf. \cite[proposition 2.1]{FriJa}, à savoir que les $H$-périodes \eqref{eq:intro-H per} sont nulles pour $p>q$ et toute fonction cuspidale. Si $n=p\geq 2$, cet énoncé résulte immédiatement de l'orthogonalité entre formes cuspidales et les séries d'Eisenstein,  donc  entre formes cuspidales et la fonction constante (qui est un résidu de séries d'Eisenstein). D'ailleurs, si $n=p\geq 2$, la variété $S$ est réduite au singleton $\{I_n\}$ et la  fonction $K_\Phi$ est constante: seule la représentation triviale peut apparaître dans sa  décomposition spectrale ce qui corrobore le  fait que le morphisme $\rho$ se réduit au morphisme $SL(2,\CC)\to\Gc$ donné par  la représentation $\Sym^{n-1}$ de $SL(2,\CC)$.

 Supposons maintenant $p=q$. Soit $(\pi,V_\pi)$ une représentation cuspidale de $G(\AAA)$ de caractère central trivial. Friedberg et Jacquet montrent le théorème suivant:

 \begin{theoreme}\cite[proposition 2.3, théorème 4.1]{FriJa}
La $H$-période sur $V_\pi$ est non nulle si et seulement si la fonction $L$ de carré extérieur $L(\pi,\La^2,s)$ a un pôle en $s=1$ et si la fonction $L$ standard $L(\pi,s)$ est non nulle en $s=1/2$.   
\end{theoreme}

 Conjecturalement, la condition sur la fonction de carré extérieur signifie que le paramètre de Langlands  de $\pi$ se factorise à travers le morphisme $\rho:\Gc_S=Sp(2p,\CC)\to GL(n,\CC)$ (ici $n=2p$) ce qui est conforme aux prédictions du  § \ref{S:intro-gp dual}. Nous dirons aussi que la représentation $\pi$ est symplectique.
 \end{paragr}

 \begin{paragr}[Cas $H=GL(p+1)\times GL(p)$ avec $p\geq 1$.] --- \label{S:(p+1,p)}Dans ce cas, on a vu que la composante cuspidale de $K_\phi(g)=\int_{[H]}K_f(g, h)\,dh$ est nulle. Soit $\Xgo(G)$ l'ensemble des données cuspidales de $G$ (on renvoie le lecteur à \cite[section 12]{Ar-cours} pour la notion de donnée cuspidale). On a alors une décomposition hilbertienne due à Langlands, cf.  \cite[lemme 12.4]{Ar-cours},
    \begin{align*}
     L^2([G])=\bigoplus_{\chi\in \Xgo(G)} L^2_\chi([G])
   \end{align*}
   avec $[G]=G(F)\back G(\AAA)$ où chaque facteur est stable par l'action de $G(\AAA)$. Pour toute donnée $\chi\in \Xgo(G)$, on note  $K_{\chi,f}$ le noyau de l'opérateur de convolution de $f$ sur le facteur correspondant à $\chi$.

Soit $P\subset G$ le sous-groupe parabolique standard de $G$ de type $(2p,1$ et soit  $M=GL(2p)\times GL(1)$ son facteur de Levi standard.
   
      \begin{theoreme}\label{thm:intro-perK} (cf. théorème \ref{thm:chi dvpt})
   Soit $\sigma$ une représentation automorphe cuspidale de $GL(2p)$ et  $\pi=\sigma\otimes 1$ qui est une représentation automorphe cuspidale de $M$. Soit  $\chi\in \Xgo(G)$ la donnée cuspidale définie par le couple $(M,\pi)$.  L'application
        \begin{align}\label{eq:intro-perK}
g\in G(\AAA)\mapsto         \int_{[H]}K_{\chi,f}(g, h)\,dh
   \end{align}
   est non identiquement nulle sur $G(\AAA)$ si et seulement si  $\sigma$ est symplectique c'est-à-dire si la fonction $s\in \CC\mapsto L(\sigma,\La^2,s)$ a un pôle en $s=1$.  De plus, on a pour tout $x\in G(\AAA)$
   \begin{align*}
        \int_{[H]}K_{\chi,f}(g, h)\,dh=\sum_{\varphi\in \bc_\pi}   E(x,I_P(f)\varphi,0) \overline{J_P(\varphi) }.
   \end{align*}
      \end{theoreme}

La somme ci-dessus porte sur une certaine base orthormale $\bc_\pi$ de la représentation induite, notée $I_P$, de la représentation $\pi$, cf. § \ref{S:car rel}. On obtient ainsi un \og caractère relatif\fg{} construit ici à l'aide d'une série d'Eisenstein  $E(x,\varphi,\la)$ évaluée en le paramètre $\la=0$ et d'une forme linéaire notée  $J_P$   sur la représentation induite: cette dernière s'appelle une période d'entrelacement et elle est définie  au § \ref{S:JP}.  On observera que,   bien que la décomposition spectrale du noyau  $K_{\chi,f}$ lui-même soit purement contine, la décomposition spectrale de  son intégrale sur $[H]$  est discrète.  Notons enfin que le théorème  \ref{thm:intro-perK} est conforme aux conjectures du § \ref{S:intro-gp dual} à savoir que, pour que l'expresssion \eqref{eq:intro-perK} soit non identiquement nulle, il faut que le paramètre de $\pi$ se factorise par $\rho$.
   
 \end{paragr}

 \subsection{Une contribution spectrale dans la formule des traces de Guo-Jacquet}\label{ssec:intro-spec}

 \begin{paragr}On reprend les notations du § \ref{S:intro-GJ}. On vient de voir des  liens  entre la $H$-distinction et les valeurs spéciales de certaines fonctions $L$ dans le cas du groupe $G=GL(n)$ pour lequel on dispose (cf. \cite{FriJa}) de la construction de ces  fonctions $L$ à l'aide d'intégrales zêta. On a $n=p+q$. Supposons de plus qu'on a $p=p'd$, $q=q'd$  et que  $D$ est une algèbre à division centrale sur $F$ de dimension  $d^2$ ; on peut comme auparavant définir un élément $\theta$  d'ordre $2$ de $G'=GL_{p'+q'}(D)$ dont le centralisateur est  le sous-groupe $H'=GL_{p'}(D)\times GL_{q'}(D)$. On peut étudier le spectre automorphe de l'espace symétrique associé  en terme de la correspondance de Jacquet-Langlands entre $G'$ et $G$ établie en toute généralité par Badulescu et Renard, cf. \cite{Badu,BaduR}.  Il est par exemple tentant de prédire des liens entre la $H'$-distinction d'une représentation automorphe $\pi'$ de $G'$ et la $H$-distinction de son  transfert de Jacquet-Langlands à $G$. On renvoie à   \cite{matringe2025} pour des résultats partiels dans cette direction lorsque $p=q$ qui reposent sur l'étude de périodes régularisées de séries d'Eisenstein. Une autre façon d'aborder cette question est d'utiliser une comparaison de formules des traces relative (comme suggéré dans \cite{Zha2}). La philosophie de cette approche est qu'il devrait exister des relations entre intégrales orbitales relatives (à l'instar du \og transfert\fg{} pour la formule des traces d'Arthur) duales de relations entre caractères relatifs construits respectivement à l'aide des $H'$-périodes et des $H$-périodes. Cela motive, on l'espère,  la construction d'une formule des traces relatives, dite de Guo-Jacquet  dans ces cas-là,   entamée dans notre travail \cite{CLi} avec Huajie Li.
    \end{paragr}

    \begin{paragr}[Développement spectral de la formule des traces de Guo-Jacquet.] ---  \label{S:intro-DSFT}On se place désormais dans la situation du § \ref{S:intro-GJ} dont on reprend les notations.  Nous allons maintenant brièvement décrire la construction de cette formule des traces dans ce cadre. Soit $f\in \Sc(G(\AAA))$ une fonction de Schwartz. Le principe de la formule des traces relatives est \emph{a priori} d'intégrer sur $Z_G^\theta(\AAA) H(F)\back H(\AAA)$ le noyau \eqref{eq:intro-KPhi} et d'obtenir deux développements, l'un selon les doubles classes dans $H(F)\back G(F)/H(F)$ et l'autre selon les données cuspidales dans $\Xgo(G)$. Comme on le sait bien, cette intégrale n'est pas convergente et on va, avant d'intégrer, remplacer le noyau $K_f$ par un noyau modifié. Celui-ci dépend d'un point $T$  dans l'espace réel $\ago_{P_0}=\Hom_\ZZ(X^*(P_0),\RR)\simeq \RR^n$ associé au groupe des caractères $X^*(P_0)$ du sous-groupe $P_0$ des matrices triangulaires supérieures dans $G$. Soit $M_0$ le sous-groupe de $P_0$ formé des matrices diagonales.  Dans la suite, on identifie le groupe de Weyl $W$ du couple $(G,M_0)$ au sous-groupe de $G(F)$ des matrices de permutation. On suppose que ce point est assez positif c'est-à-dire assez profond dans la chambre de Weyl aiguë associée à $P_0$. On définit le noyau modifié $K^{T}_{\chi,f}$, associé à la donnée cuspidale $\chi\in \Xgo(G)$,  par la formule
      \begin{align}\label{eq:introKT}
        K^{T}_{\chi,f}(x,y)=\sum_{P_0\subset P\subset G} \eps_P^G \sum_{w_1 \in   \, _PW_{H}} \sum_{\delta_1\in P_{w_1}^{H}(F) \back H(F)}  \hat\tau_P(H_P(w_1\delta_1x)-T) \times \\ \nonumber \left[\sum_{w_2\in   \, _PW_H} \sum_{\delta_2\in P_{w_2}^H(F) \back H(F)} K_{P,\chi,f}(w_1\delta_1x,w_2\delta_2y)\right]
      \end{align}
      pour tout $x,y\in G(\AAA)$. Expliquons les notations: la première somme porte sur les sous-groupes paraboliques $P$ qui contiennent $P_0$. À ceux-ci on associe:
      \begin{itemize}
      \item un $\RR$-espace vectoriel $\ago_P=\Hom(X^*(P),\RR)$ comme ci-dessus; on a d'ailleurs une décomposition naturelle $\ago_P=\ago_P^G\oplus \ago_G$ où $\ago_P^G$ est le noyau de la projection canonique $\ago_P\to \ago_G$;
      \item  $\eps_P^G =(-1)^{\dim(\ago_P^G)}$;
      \item une décomposition de Levi $P=M_P N_P$ où $N_P$ est le radical unipotent et $M_P$ l'unique  facteur de Levi de $P$ qui contient $M_0$ ;
      \item le sous-ensemble $\, _PW_{H}\subset W$ formé des $w\in W$ tels que  $M_P\cap P_0=M_P\cap wP_0 w^{-1}$  et $P_0\cap H\subset w^{-1}Pw$ ;
      \item pour $w\in  \, _PW_{H}\subset W$, le sous-groupe parabolique de $H$ défini par $P_{w}^{H}=w^{-1}Pw\cap H$;
      \item $\hat \tau_P$,  la fonction sur $\ago_P$, caractéristique de la chambre ouverte obtuse associée à $P$ (\og duale des poids \fg{});
      \item $H_P:G(\AAA)\to \ago_P$, l'application,  invariante à droite par $K$ le sous-compact maximal usuel de $G(\AAA)$,  qui coïncide sur $P(\AAA)$ avec la composition de l'application canonique $P(\AAA)\times X^*(P)\to \AAA^\times$ et du logarithme du  module usuel;
      \item $K_{P,\chi,f}$, le noyau de l'opérateur donné par l'action par convolution de $f$ sur le facteur $L^2_\chi(M_P(F)N_P(\AAA)\back G(\AAA))\subset L^2(M_P(F)N_P(\AAA)\back G(\AAA))$ (cf. \cite[section 12]{Ar-cours} et  \cite[§ 2.9.2]{BCZ}).
      \end{itemize}

      Toutes les sommes qui apparaissent dans     \eqref{eq:introKT}  sont finies (ou à support fini qui dépend de $x$) excepté la  somme intérieure entre crochets qui est néanmoins convergente, cf. \cite[remarque 3.2.2.1]{CLi}. Notons aussi que le noyau modifié ne dépend que de la projection $T^G$ de $T$ sur $\ago_P^G$.      Un premier résultat, \cite[théorème 3.2.3.2]{CLi} est que l'expression
      \begin{align*}
        \sum_{\chi\in \Xgo(G)} | K^{T}_{\chi,f}(x,y)|
      \end{align*}
      définit une fonction sur $H(\AAA)\times H(\AAA)$ invariante par le plongement diagonal du centre $Z_G(\AAA)$  et qui est à décroissance rapide sur le quotient. On peut donc définir une forme linéaire sur $\Sc(G(\AAA))$ par
      \begin{align*}
        J_\chi^T(f)=\int_{[H]_G} \int_{[H]}  K^{T}_{\chi,f}(x,y)\, dx dy
      \end{align*}
      où
      \begin{align*}
              [H]_G=Z_G(\AAA)H(F)\back H(\AAA).
      \end{align*}
      En fait, cette forme linéaire est continue pour la topologie naturelle sur $\Sc(G(\AAA))$, cf. \cite[théorème 3.3.2.1]{CLi}. 
      À ce stade, on peut analyser la dépendance en $T$. Il se trouve que l'application $T\mapsto J_\chi^T(f)$ coïncide dans un cône avec une exponentielle-polynôme en $T$; on  définit alors $J_\chi(f)$ comme le terme constant de cette exponentielle-polynôme, cf. \cite[proposition 3.5.3.2]{CLi}. L'espace de Schwartz est muni d'actions par translations à droite et à gauche de $G(\AAA)$ donc, par restriction de $H(\AAA)$. Pour l'action à droite, la distribution $f\mapsto J_\chi(f)$ est $H(\AAA)$-invariante. Il n'en est pas de même, en général, pour l'action à gauche, cf. \cite[proposition 3.5.5.1]{CLi}. C'est analogue à ce qui se passe pour la formule des traces d'Arthur, qui donne naissance \emph{a priori} à des distributions non-invariantes. La somme
      \begin{align*}
        \sum_{\chi\in \Xgo(G)} J_\chi(f)
      \end{align*}
      converge absolument. Cela donne  le développement spectral de la formule des traces de Guo-Jacquet selon les données cuspidales. Il y a aussi un développement géométrique établi dans \cite[section 4]{CLi} qui ne nous concernera pas ici.
    \end{paragr}

    \begin{paragr}[Opérateur de troncature.] ---\label{S:intro-LaTtheta}
      Une question qui se pose ensuite est d'obtenir la décomposition spectrale de chaque distribution $J_\chi(f)$. Pour cette question, il est plus facile de remplacer le noyau modifié par le noyau initial auquel on applique un opérateur de troncature. Depuis son introduction par Arthur dans \cite{ar2} et l'article de Jacquet-Lapid-Rogawski \cite{JLR},  de nombreux auteurs ont proposé des variantes  de l'opérateur de troncature d'Arthur. Celui que nous introduisons et notons $\La^T_\theta$  dépend du paramètre $T $ déjà utilisé ci-dessus et il est donné par la formule suivante: pour toute fonction $\varphi$, disons continue, sur $[G]$ et tout $x\in G(\AAA)$, 
       \begin{align}\label{eq:intro-LaTtheta}
    (\Lambda^{T}_\theta\varphi)(x)= \sum_{P_0\subset P\subset G}  \eps_P^G \sum_{w \in  \, _PW^G_H}   \sum_{\delta \in P_{w}^H(F)\back H(F)}\hat\tau_P(H_P(w\delta x)-T)  \varphi_P(w\delta x),
  \end{align}
où l'on utilise le terme constant le long de $P$ défini par
\begin{align}\label{eq:intro-terme cst}
\forall g\in G(\AAA), \ \  \varphi_P(g)=\int_{[N_P]} \varphi(ng)\, dn
\end{align}
avec $dn$  la mesure invariante sur $[N_P]=N_P(F)\back N_P(\AAA)$ de masse totale $1$. Les autres notations sont celles qui interviennent dans la définition  \eqref{eq:introKT} du noyau modifié. 

\begin{remarque}\label{rq:op-zy}
  Même s'il est formellement proche de l'opérateur que l'on peut déduire des constructions générales de Zydor dans   \cite[section 3.7]{Zydor}, il est important de souligner que notre opérateur n'est pas celui de Zydor. En particulier, notre opérateur vérifie que  $\Lambda^{T}_\theta\varphi$ converge simplement vers $\varphi$ lorsque $T$ tend vers l'infini dans la chambre positive. Celui de Zydor ne vérifie pas cette propriété; on pourra consulter  \cite[remarque 2.4.2.1]{CLi}  pour plus d'explication.
\end{remarque}

Pour une description des principales propriétés de l'opérateur   $\Lambda^{T}_\theta$, on renvoie à \cite[§ 2.4.3]{CLi}. Indiquons cependant qu'il transforme les fonctions à croissance uniformément modérée sur $[G]$ et invariantes par le centre $Z_G(\AAA)$ en des fonctions à décroissance rapide sur $Z_G(\AAA)H(F)\back H(\AAA)$. Une autre propriété cruciale est la suivante:

\begin{theoreme} \cite[théorème 3.3.3.1]{CLi}
  \label{thm:intro-asym}
  Soit $\|\cdot\|$ une norme sur $\ago_{P_0}^G$. Il existe une semi-norme  continue $\|\cdot\|_{\Sc} $ sur $\Sc(G(\AAA))$ telle que pour tout $T\in \ago_0$ suffisamment positif,  tout  $f\in \Sc(G(\AAA))$ et tout $\chi\in \Xgo(G)$, on a
    \begin{align}
&\label{eq:intro-asym} \left| J_\chi^T(f)-  \int_{[H]_G} \int_{[H]}  (\La^T_{\theta} K_{\chi,f})(x,y)\, dx dy \right|\leq e^{-\|T^G\|}  \|f\|_{\Sc}.
    \end{align}
  \end{theoreme}
  
         Ici $T$ suffisamment positif signifie que $T$ est assez profond et assez loin des murs dans la chambre de Weyl, pour une définition précise on renvoie à  \cite[§ 2.1.14]{CLi}. La notation $\La^T_{\theta} K_{\chi,f}$ indique qu'on a appliqué l'opérateur $\La^T_\theta$ à la variable de gauche du noyau $K_{\chi,f}$.  Des majorations du noyau $K_{\chi,f}$ issues de \cite[lemme 2.10.1.1]{BCZ} et les propriétés de l'opérateur  $\La^T_{\theta}$, cf. \cite[proposition 2.4.3.1]{CLi}, assurent que l'intégrale de $\La^T_{\theta} K_{\chi,f}$ dans \eqref{eq:intro-asym} converge absolument. Le théorème ci-dessous est précieux pour obtenir un développement spectral fin de $J_\chi(f)$ et ainsi une expression plus explicite. Dans cette note, on illustre cette approche dans un cas particulier instructif.
\end{paragr}

\begin{paragr}[Calcul d'une contribution $J_\chi(f)$.] --- On se place désormais dans la situation du § \ref{S:(p+1,p)}, c'est-à-dire on prend $G=GL(n)$ avec  $n=2p+1$ et $H=GL(p+1)\times GL(p)$.  Soit $P\subset G$ le sous-groupe parabolique  standard de $G$ de type $(2p,1)$ et soit  $M$ son facteur de Levi standard. Soit  $\pi=\sigma\otimes 1$ une représentation automorphe cuspidale de $M$ et $\chi$  la donnée cuspidale définie par $(M,\pi)$, cf. théorème \ref{thm:intro-perK}.  Pour toute fonction test $f\in \Sc(G(\AAA))$, on introduit alors le \og{}caractère relatif pondéré\fg
 \begin{align}\label{eq:intro-car pondere}
    J_{P,\pi}(f)= \sum_{\varphi\in \bc_\pi}   J_P(\mc_{P,\pi} I_P(f)\varphi) \overline{J_P(\varphi)}.
\end{align} 
La somme porte sur une base orthonormée de l'induite $I_P$ de $\pi$, cf. § \ref{S:car rel}. La forme linéaire $J_P$, une période d'entrelacement,  intervient dans le théorème  \ref{thm:intro-perK} et est définie § \ref{S:JP}. L'opérateur $\mc_{P,\pi}$ est une sorte de dérivée logarithmique d'opérateurs d'entrelacement qui intervient dans la définition des caractères pondérés d'Arthur, pour une définition précise cf. \eqref{eq:formule-mc}. Notre principal résultat est alors le suivant:

  \begin{theoreme}\label{thm:intro-Jchi}(cf. théorème \ref{thm:Jchi})
Pour tout $f\in \Sc(G(\AAA))$, on a
    \begin{align*}
      J_\chi(f)=   J_{P,\pi}(f).
    \end{align*}
  \end{theoreme}

  Le théorème \ref{thm:intro-Jchi} repose sur les théorèmes \ref{thm:intro-asym} et   \ref{thm:intro-perK} ainsi que sur l'étude des périodes de certaines  séries d'Eisenstein tronquées par l'opérateur $\La^T_\theta$, cf. section \ref{ssec:period-Eis}. Un des intérêts du théorème est qu'il offre la possibilité,  moyennant une factorisation des périodes d'entrelacement, de relier la distribution $J_\chi$ à des distributions locales, de même que les caractères pondérés d'Arthur sont reliés à des analogues locaux.
  
\end{paragr}

\subsection{Organisation de l'article}

\begin{paragr} La section  \ref{ssec:notations}   introduit les principales notations et une description d'ensembles de doubles classes.  On définit  certaines périodes d'entrelacement dans la section \ref{ssec:entrelacement}. L'étude des périodes de certaines séries d'Eisenstein tronquées apparaît à la section \ref{ssec:period-Eis}. Leur calcul, en terme de périodes d'entrelacement, est obtenu dans le théorème \ref{thm:LaTEis}  pour un paramètre complexe non nul avec, comme corollaires \ref{cor:LaTEis} et \ref{cor:LaTEisven 0}, d'une part une équation fonctionnelle pour les périodes d'entrelacement et d'autre part  le calcul pour le paramètre nul. Le reste de la section \ref{ssec:period-Eis}, à partir du § \ref{S:preuve LaTEis} est consacrée à la démonstration du théorème \ref{thm:LaTEis}. Dans la section \ref{ssec:unecontrib}, on rappelle la notion de caractère relatif et on définit le caractère relatif pondéré  $J_{P,\pi}$ et on énonce l'égalité $J_\chi=J_{P,\pi}$. On vérifie aussi certaines propriétés de covariance de ces deux distributions. La démonstration du théorème \ref{thm:intro-perK} est donnée à la section \ref{ssec:undvptspec}, celle du théorème \ref{thm:intro-Jchi} à la section finale \ref{ssec:demo-thmpal}.
  
\end{paragr}

\begin{paragr}[Remerciements.] --- Je remercie les organisateurs de la conférence \emph{Trace Formula, Endoscopic Classification and Beyond: the Mathematical Legacy of James Arthur} au Fields Institute à Toronto de m'avoir permis d'exposer une partie des résultats présentés ici.  Je remercie Huajie Li pour notre collaboration et nos nombreuses discussions sur la formule des traces de Guo-Jacquet.  Je remercie également Linli Shi de l'intérêt qu'il a manifesté pour cet article. Je remercie enfin l'Institut Universitaire de France (IUF) pour m'avoir offert d'excellentes conditions de travail durant la rédaction  de cet article. 
  \end{paragr}

\section{Un développement spectral}

\subsection{Notations}\label{ssec:notations}

\begin{paragr}
  Dans toute la suite, $F$ est un corps de nombres. Soit $V_F$ l'ensemble des places de  $F$ et $V_{F,\infty}$ et $V_F^\infty$ les sous-ensembles des places archimédiennes respectivement. Soit $\AAA$ l'anneau des adèles de $F$. On  note $|\cdot|$ la fonction module sur le groupe multiplicatif $\AAA^\times$. On espère qu'il n'y aura pas de confusion avec le module complexe sur $\CC$ également noté $|\cdot|$. On note $\Re$ et $\Im$ les parties réelles et imaginaires d'un nombre complexe. Soit $\AAA^1$ le noyau du module $|\cdot|$. On munit  $\AAA^1$ de la mesure de Haar qui donne le volume $1$ au quotient $F^\times\back \AAA^1$ muni de la mesure quotient, $F^\times$ étant muni de la mesure de comptage. Le module identifie le quotient $\AAA^1\back \AAA^\times$ à $\RR_+^\times$. On munit alors  $\AAA^\times$ de la mesure de Haar qui donne au quotient $\AAA^1\back \AAA^\times$ la mesure $dt/t$ sur $\RR_+^\times$ où $dt$ est la mesure de Lebesgue sur $\RR$.
\end{paragr}

\begin{paragr}
  Pour tout entier $n\geq 1$, on note par  $G_n$ le groupe $GL(n)$ sur le corps $F$. On note $[G_n]=G_n(F)\back G_n(\AAA)$ et cette notation vaut pour tous les groupes algébriques définis sur $F$ rencontrés. Soit $Z_n\subset T_n\subset B_n\subset G_n$  respectivement le centre, le sous-tore maximal standard et le  sous-groupe de Borel standard. Ce dernier  est, par définition, le stabilisateur dans $G_n$ du drapeau complet $\vect(e_1)\subset \vect(e_1,e_2)\subset \ldots  \subset \vect(e_1,\ldots,e_n)$ où $(e_i)_{1\leq i\leq n}$ est la base canonique de $F^n$. Soit $N_n$ le radical unipotent de $B_n$. Le tore $T_n$ est le stabilisateur dans $G_n$ des droites engendrées par les vecteurs de la base canonique. Soit $W_n$ le groupe de Weyl de $(G_n,T_n)$ ; par définition, c'est le quotient du normalisateur dans $G_n(F)$ de $T_n$ par le sous-groupe $T_n(F)$. Dans cet article, on identifie $W_n$ au sous-groupe de $G_n(F)$ des matrices de permutation (c'est-à-dire des automorphismes de $F^n$ qui stabilisent la base canonique).

  Un sous-groupe parabolique de $G_n$ est dit standard, resp. semi-standard, s'il contient $B_n$, resp. $T_n$. Pour tout $1\leq k \leq n-1$, soit  $P_{k,n-k}$ le sous-groupe parabolique standard de $G_n$ de type $(k,n-k)$ c'est-à-dire le stabilisateur du sous-espace $\vect(e_1,\ldots,e_k)$.  Soit $N_{k,n-k}$ le radical unipotent de $P_{k,n-k}$.  Soit $\bar P_{k,n-k}$ le sous-groupe parabolique opposé à $P_{k,n-k}$: par définition c'est le stabilisateur du sous-espace $\vect(e_{k+1},\ldots,e_n)$. 
  Pour tout sous-groupe parabolique semi-standard $P$ de $G_n$, on note $M_{P}$ l'unique facteur de Levi qui contient $T_n$ et $N_P$ le radical unipotent de $P$. On a donc une décomposition de Levi $P=M_P N_P$. Soit $Z_{M_P}\subset M_P$ le centre de $M_P$.

  On note $K_n=\prod_{v\in V_F} K_{n,v}$ le sous-groupe compact maximal de $G_n(\AAA)$ défini ainsi: si $v$ est réelle,  $K_{n,v}\subset G_n(\RR)$ est le groupe orthogonal pour la forme quadratique sur $\RR^n$ de base orthonormale  $(e_i)_{1\leq i\leq n}$;  si $v$ est complexe,  $K_{n,v}\subset G_n(\CC)$ est le groupe unitaire pour la forme hermitienne  sur $\CC^n$  de base orthonormale    $(e_i)_{1\leq i\leq n}$; enfin si $v$ est non-archimédienne, $K_{n,v}\subset G_n(F_v)$ est le sous-groupe des automorphismes du réseau $\oplus_{1\leq i\leq n} \oc_v e_i$.
\end{paragr}

\begin{paragr}[Mesures de Haar.] ---  Les groupes $Z_n(\AAA)\simeq \AAA^\times$  et  $T_n(\AAA)\simeq (\AAA^\times)^n$ héritent de la mesure de Haar sur $\AAA^\times$. Pour tout groupe unipotent $N$ défini  sur $F$, on munit $N(\AAA)$ de la mesure de Haar qui donne le volume $1$ à $[N]$. On munit $K_n$ de la mesure de Haar qui donne $1$ comme volume total $1$. On munit alors $G_n(\AAA)$ de la mesure de Haar telle que, pour toute fonction continue à support compact, on ait:
  \begin{align*}
    \int_{G_n(\AAA)} f(g)\,dg=\int_{T_n(\AAA)}\int_{N_{B_n}(\AAA)} \int_{K_n} f(tnk)\,dtdndk.
  \end{align*}

  Notons que les sous-groupes de Levi standard de $G_n$, qui sont des produits de groupes $G_r$ pour $r\leq n$,  ainsi que leur centre sont alors  aussi munis de mesures de Haar. Soit $H\subset G_n$ un sous-groupe qui contient le centre $Z_n$ de $G_n$. On pose
  \begin{align*}
    [H]_{G_n}= Z_n(\AAA) H(F)\back H(\AAA).
  \end{align*}
 La mesure de Haar sur le groupe $Z_n(\AAA) H(F)$ est obtenue à l'aide de l'isomorphisme     $Z_n(\AAA) H(F)\simeq (Z_n(\AAA)\times H(F))/Z_n(F)$ et de la mesure quotient sur $(Z_n(\AAA)\times H(F))/Z_n(F)$, les groupes $H(F)$ et $Z_n(F)$ étant munis de la mesure de comptage. Alors, si $H(\AAA)$ est muni explicitement  d'une mesure de Haar, on munit  $[H]_{G_n}$ de la mesure quotient.
\end{paragr}

\begin{paragr} On fixe un entier $p\geq 1$ et  $n=2p+1$. Soit $G=G_n$ , $P_0=B_n$ et $M_0=T_n$ et $W=W_n$ identifié au sous-groupe des matrices de permutation. Pour toute place $v\in V_F$, soit  $K_v=K_{n,v}$ et soit $K=\prod_{v\in V_F} K_v$. 

  Soit $H=G_{p+1}\times G_p$ le sous-groupe de $G$ qui stabilise les sous-espaces $\vect(e_1,\ldots,e_{p+1})$ et $\vect(e_{p+2},\ldots,e_{2p+1})$. Soit $P$ et $Q$ les sous-groupes paraboliques standard  de $G$ de type respectif $(1,2p)$ et $(2p,1)$.    Soit $M$ et $L$ les facteurs de Levi standard respectivement de $P$ et $Q$. Ainsi $M= G_1\times G_{2p}$ est le  sous-groupe de $G$ qui stabilise les sous-espaces $\vect(e_1)$ et $\vect(e_{2},\ldots,e_{2p+1})$ alors que $L=G_{2p}\times G_1$ est celui qui stabilise les sous-espaces $\vect(e_{2p+1})$ et $\vect(e_{1},\ldots,e_{2p})$.
\end{paragr}

\begin{paragr}
  On fixe une hauteur $\|\cdot\|:Z_G(\AAA)\back G(\AAA) \to \RR^\times_+$, cf. \cite[I.2.2]{MWlivre}. On pose $\|g\|_G=\inf_{\delta\in G(F)}  \|\delta g\|$ pour tout $g\in [G]_G=Z_G(\AAA)G(F)\back G(\AAA)$. Plus généralement, pour tout $k\geq 1$, on définit de la même façon $\|\cdot\|_{G_k}$ sur $[G_k]_{G_k}$.  On pose  aussi $\|g\|_H=\inf_{\delta\in H(F)}  \|\delta g\|$ pour tout $g\in [H]_G=Z_G(\AAA)H(F)\back H(\AAA)$.  Il est clair que pour tout $g\in [H]_G$ on a $\|g\|_G\leq \|g\|_H$. En fait, il existe $c,N>$ tel que pour tout $g\in [H]_G$ on a $\|g\|_H\leq c\|g\|_G^N$, \cite[Proposition A.1.1]{RBP}.  Pour des majorations sur  $[H]_G$ pour lesquelles l'exposant de $\|\cdot\|_G$ ou $\|\cdot\|_H$ importe peu, on utilisera indifféremment   $\|\cdot\|_G$ ou $\|\cdot\|_H$. 
\end{paragr}

\begin{paragr}   Comme au § \ref{S:intro-DSFT}, on attache à $P$ un espace vectoriel réel $\ago_P=\Hom(X^*(P),\RR)$. On dispose d'une application $H_P:G(\AAA)\to \ago_P$ qui induit un morphisme surjectif $[Z_M]\to \ago_P$. Le groupe $[Z_M]$ est muni de la mesure quotient. Le noyau de ce morphisme est muni de la mesure de Haar qui donne un volume total $1$. On munit alors $\ago_P$ de la mesure quotient. De même, on munit $\ago_G$ d'une mesure de Haar. et $\ago_P^G\simeq \ago_P/\ago_G$ est muni de la mesure quotient.

  Soit $\ago_P^*=X^*(P)\otimes \RR$. On a une dualité naturelle $\bg \cdot,\cdot\bd: \ago_P^* \times \ago_P \to \RR$. Dualement à la décomposition  $\ago_P^G\oplus \ago_G=\ago_P$ on a une décomposition   $\ago_P^{G,*}\oplus \ago_G^*=\ago_P^*$ où  $\ago_G^*=X^*(G)\otimes \RR$. Le $\RR$-espace vectoriel $i\ago_{P}^{G,*}$ est alors muni de la mesure de Haar duale au sens où pour toute fonction $\phi$ dans la classe de Schwartz de $\ago_P^G$ on a
  \begin{align*}
\int_{i\ago_{P}^{G,*}} \int_{\ago_{P}^G}  \phi(H) \exp (-\bg \la ,H\bd)\,dHd\la=\phi(0) .
  \end{align*}
  Soit $\al\in \ago_P^{G,*}$ l'unique élément tel que $z\in Z_M(\AAA)$ agisse sur $N_P(\AAA)$ par le caractère $z\mapsto \exp(\bg \al, H_P(z)\bd)$. On utilise aussi $\rho_P^G$ défini par $2\rho_P^G=2p \al$. Soit $\al^\vee\in \ago_P^{G}$: c'est l'image par l'application canonique $\ago_{P_0}\to \ago_P$ de l'unique coracine simple de $M_0$ dans $N_P$.  Pour tout $\la\in \ago_P^{G,*}$, on pose
 \begin{align}\label{eq:theta_P}
    \theta_P(\la)= \vol(\ago_P^G / \ZZ \al^\vee)^{-1} \bg \la, \al^\vee\bd.
 \end{align}

  Toutes ces notations valent pour le sous-groupe parabolique $Q$; l'analogue  de $\al$  pour $Q$ est noté $\be$ et on a $2\rho_Q^G=2p \be$.

\end{paragr}

\begin{paragr}[Doubles classes.] --- On rappelle qu'on a défini, au § \ref{S:intro-DSFT}, des  sous-ensembles  $  _P W_H$ et $  _Q W_H$ de $W$. Pour $w\in _P W_H$, on pose $P_w=w^{-1}Pw$ et $P_w^H=P_w\cap H$. Pour $w$ l'élément trivial, on omet l'indice $w$ de sorte qu'on a $P^H=P\cap H$. Toutes ces notations valent aussi pour $Q$.

  \begin{lemme}\label{lem:classes}
    \begin{enumerate}
    \item     L'ensemble  $_P W_H$ est formé de deux éléments à savoir l'élément trivial de $W$ et l'élément
    \begin{align*}
      w_1=
      \begin{pmatrix}
        0 & 1 & 0 \\ I_{p+1} & 0 & 0\\ 0 & 0 &I_{p-1}
      \end{pmatrix}.
    \end{align*}
    On a $P^H=P_{1,p} \times G_p$ et $P_{w_1}^H= G_{p+1}\times P_{1,p-1}$.
  \item Le quotient $P(F)\back G(F)/ H(F)$ a trois éléments: un système de représentants est donné par l'ensemble  $_P W_H$ auquel on adjoint l'élément
    \begin{align*}
       w_2=
      \begin{pmatrix}
        0 & 0 & 1 \\ 1 & 0 & 1\\ 0 & I_{2p-1} &0
      \end{pmatrix}.
    \end{align*}
  \item   L'ensemble  $_Q W_H$ est formé de deux éléments à savoir l'élément trivial de $W$ et l'élément
    \begin{align*}
      w_3=
      \begin{pmatrix}
        I_{p}  & 0 & 0 \\ 0 & 0 &   I_{p}\\ 0 & 1 & 0
      \end{pmatrix}.
    \end{align*}
    On a $ Q^H=G_{p+1} \times P_{p-1,1}$ et $Q_{w_3}^H= P_{p,1}\times G_p$.
  \end{enumerate}
\end{lemme}

\begin{preuve}
  1. Soit $w\in   \,_PW_H$. Alors $P_w$ est le stabilisateur de la droite engendrée par $w^{-1}e_1$. Pour que $H\cap P_w$ contienne $H\cap P_0$,  il faut et il suffit que  $w^{-1}e_1\in \{e_1,e_{p+2}\}$. Supposons d'abord $we_1=e_1$. On a donc $w\in W\cap M(F)$; il s'ensuit que $w$ est trivial puisqu'on a $M\cap (w P_0 w^{-1})=M\cap P_0$.    Supposons ensuite $w^{-1} e_1=e_{p+2}$.  Comme $M\cap (w P_0 w^{-1})=M\cap P_0$, on voit alors que $wP_0w^{-1}$ est le stabilisateur du drapeau complet
   \begin{align*}
     \vect(e_2)\subset \vect(e_2,e_3) \subset \ldots  \subset\vect(e_2,\ldots, e_{p+2})  \subset  \vect(e_1,e_2,\ldots, e_{p+2}) \subset  \vect(e_1,\ldots, e_{p+3}) \subset \ldots 
   \end{align*}
   Par conséquent, $w=w_1$. Le reste des assertions est évident.

   2. C'est essentiellement \cite[lemme 5.1]{FriJa}. Notons que l'élément $\xi_1$ de \cite[lemme 5.1]{FriJa} appartient à la classe $w_1 (W\cap H(F))$.
   
   3. se démontre comme 1.
 \end{preuve}

\end{paragr}

\subsection{Périodes d'entrelacement}\label{ssec:entrelacement}

\begin{paragr}
  Soit $\sigma$ une représentation automorphe cuspidale irréductible de $G_{2p}(\AAA)$ de caractère central trivial. Soit $\pi$ la  représentation automorphe cuspidale de $M(\AAA)$ donnée par le produit tensoriel externe de la représentation triviale de $G_1(\AAA)$ et $\sigma$.  Soit $\pi'$ la  représentation automorphe cuspidale de $L(\AAA)$ donnée par le produit tensoriel externe  de $\sigma$ et de la représentation triviale de $G_1(\AAA)$.
\end{paragr}

\begin{paragr} Dans ce §, on  suit \cite[section 2.6]{chaudouardsymmetric} pour les notions de formes automorphes qu'on utilise.  La représentation $\pi$ détermine un espace $\Ac_{\pi}(M)$ de formes automorphes sur lequel $M(\AAA)$ agit par translation à droite. Ces formes sont des fonctions complexes lisses sur $[M]$, invariantes sous le centre $Z_M(\AAA)$, à croissance modérée et finies sous l'action du centre de  l'algèbre enveloppante de l'algèbre de Lie  complexifiée de $M$. En revanche, on ne les suppose pas finies sous l'action d'un sous-groupe compact maximal fixé de $\prod_{v\in V_{F,\infty}} M(F_v)$. On introduit ensuite un espace $\Ac_{P,\pi}(G)$ de fonctions lisses $\varphi$ sur le quotient
  \begin{align*}
    [G]_P=Z_M(\AAA)M(F) N_P(\AAA)\back G(\AAA),
  \end{align*}
à croissance modérée, finies sous l'action du centre de  l'algèbre enveloppante de l'algèbre de Lie  complexifiée de $G$ telles que pour tout $g\in G(\AAA)$ l'application $m\in[M]\mapsto  \exp(-\bg\rho_P^G,H_P(m)\bd)  \varphi(mg)$  appartient à $\Ac_{\pi}(M)$. Cet espace  $\Ac_{P,\pi}(G)$ est muni d'une topologie, cf. \cite[§  2.6.13]{chaudouardsymmetric}.

      Soit $\varphi\in \Ac_{P,\pi}(G)$. Pour tout $\la\in \ago_{P,\CC}^{G,*}$, soit $\varphi_\la$ défini par
  \begin{align*}
    \forall g\in G(\AAA), \ \     \varphi_\la(g)=\exp(\bg\la,H_P(g)\bd)  \varphi(g).
  \end{align*}
  L'espace vectoriel  $\Ac_{P,\pi}(G)$ est muni de la norme hilbertienne de Petersson donnée par
  \begin{align*}
    \|\varphi\|_P^2=\int_{[G]_P}    |\varphi(g)|^2\, dg.
  \end{align*}

  \begin{remarque}\label{rq:mes inv}
    Dans l'intégrale ci-dessus, le groupe   $Z_M(\AAA)M(F) N_P(\AAA)$ est muni de la mesure de Haar à droite donnée par
    \begin{align*}
      \phi\mapsto \int_{N_P(\AAA)} \int_{[Z_M]} \sum_{\ga \in M(F)}  f(n z \ga)\, dzdn.
    \end{align*}
 Cette mesure de Haar n'est cependant pas une mesure de Haar à gauche: le caractère modulaire est donné par $\delta_P(x)= \exp(\bg 2\rho_P^G, H_P(x)\bd)$. Comme il est d'usage, la notation intégrale désigne en fait une forme linéaire invariante à droite sur l'espace des fonctions continues sur $G(\AAA)$ et $\delta_P$-équivariante à gauche  sous l'action de $Z_M(\AAA)M(F) N_P(\AAA)$.
  \end{remarque}
  
  Pour tout $\la\in \ago_{P,\CC}^{G,*}$, on définit une action à gauche  $I_{P,\pi}(\la)$ de $G(\AAA)$ sur $ \Ac_{P,\pi}(G)$ de la façon suivante:
  \begin{align*}
    \forall \varphi \in \Ac_{P,\pi}(G), \ \forall x\in G(\AAA), \ \ I_P(\la,x)\varphi =   (\varphi_\la(\cdot x))_{-\la}.
  \end{align*}
 L'intégration sur $G(\AAA)$ contre une fonction $f\in \Sc(G(\AAA))$ définit l'opérateur $I_P(\la,f)$.
  Pour $\la=0$, on pose $I_P=I_P(0)$ et on retrouve l'action par translation à droite de $G(\AAA)$ sur $ \Ac_{P,\pi}(G)$. Pour $\la\in i\ago_{P}^{G,*}$ cette action est unitaire pour la norme de Petersson.

  On utilisera également sans plus de commentaire les constructions analogues relatives à $Q$ et $\pi'$ avec des notations qu'on espère évidentes.
\end{paragr}

\begin{paragr}[Période d'entrelacement $J_P$.] --- \label{S:JP}La  décomposition de Levi $P=M N_P$  induit une décomposition de Levi $P^H=M^HN_P^H$ de facteur de Levi $M^H=H\cap M$ et de radical unipotent $N_P^H=H\cap N_P$ pour le  sous-groupe parabolique $P^H$ de $H$. Le  caractère $x\in P(\AAA)\mapsto \exp(\bg \rho_P^G,H_P(x)\bd$ et le  caractère modulaire $\delta_{P^H}$ du  groupe $P^H(\AAA)$  ont la même restriction à $Z_M(\AAA)M^H(F)N_P^H(\AAA)$. Pour tout $\varphi\in \Ac_{P,\pi}(G)$ on introduit  alors la période d'entrelacement
  \begin{align}\label{eq:periode entrelac}
    J_P(\varphi)=\int_{ [H]_P}  \varphi(g)\, dg,
  \end{align}
  où l'on pose  $[H]_P= Z_M(\AAA)M^H(F)N_P^H(\AAA) \back H(\AAA) $ et   l'intégrale doit être comprise dans le sens de la remarque \ref{rq:mes inv}. La proposition suivante est essentiellement une reformulation de résultats de Friedberg-Jacquet.

\begin{proposition}\label{prop:periode entrelac}(Friedberg-Jacquet)
    \begin{enumerate}
    \item L'intégrale dans le membre de droite de \ref{eq:periode entrelac} est absolument convergente et définit une forme linéaire, notée $J_P$, sur  $\Ac_{P,\pi}(G)$ qui est  $H(\AAA)$-invariante  et continue.
    \item  La forme linéaire  $J_P$ est non identiquement nulle sur $\Ac_{P,\pi}(G)$ si et seulement si la représentation $\sigma$ est symplectique.
    \end{enumerate}
  \end{proposition}

  \begin{preuve}
   On écrit $g\in H(\AAA)$ comme un couple $(g',g_3)\in (G_{p+1}\times G_p)(\AAA)$. On utilise alors la décomposition d'Iwasawa de $g'$ relativement au sous-groupe parabolique $P_{1,p}$ ce qui donne:
    \begin{align*}
      g'=  \begin{pmatrix}
        1 & u \\ 0 & I_p
      \end{pmatrix}
      \begin{pmatrix}
        g_1 & 0 \\ 0 &g_2
      \end{pmatrix} k
    \end{align*}
    avec $u\in \AAA^{p} $ (identifié à une matrice ligne), $k\in K_{p+1}$, $g_1\in G_1(\AAA)$ et $g_2\in G_p(\AAA)$. On obtient alors la formule:
    \begin{align*}
    \exp( \bg \rho_{P}^{G}, H_{P}(g)\bd)  \delta_{P^H}(\begin{pmatrix}
        g_1 & 0 & 0\\ 0 &g_2 &0 \\ 0 & 0 & g_3
      \end{pmatrix})^{-1}=  |\det(g_2) \det(g_3)^{-1}|^{\frac12 }.
    \end{align*}
    Il s'ensuit qu'on a
    \begin{align*}
      J_P(\varphi)&= \int_{ [G_p\times G_p]_{G_{2p}}} \int_{K_{p+1}}\varphi_{-\rho_P^G}\left(     \begin{pmatrix}
          1 & 0 &0\\ 0 &g_2& 0 \\ 0 & 0 & g_3
      \end{pmatrix}   \begin{pmatrix}
          k & 0 \\ 0 &I_p
        \end{pmatrix} \right)  \left|\frac{\det(g_2)}{\det(g_3)}\right|^{\frac12}dk dg_1dg_2dg_3.
         \end{align*}
         L'intégrale ci-dessus apparaît dans le travail de Friedberg-Jacquet,  \cite[théorème 5.1]{FriJa}. Il résulte de la cuspidalité de $\pi$ que, pour tout $N>0$, on peut trouver une semi-norme continue $\|\cdot \|_N$ sur $\Ac_{P,\pi}(G)$ telle que pour tout $(g_1,g_2)\in [G_p\times G_p]_{G_{2p}}$, tout $\varphi\in \Ac_{P,\pi}(G)$  et tout $k\in K_{p+1}$ on a:
         \begin{align*}
           \left| \varphi_{-\rho_P^G}\left(     \begin{pmatrix}
          1 & 0 &0\\ 0 &g_2& 0 \\ 0 & 0 & g_3
      \end{pmatrix}   \begin{pmatrix}
          k & 0 \\ 0 &I_p
        \end{pmatrix} \right)\right| \leq  \|(g_2,g_3)\|_{G_{2p}}^{-N} \|\varphi\|_N.
         \end{align*}
La convergence et la continuité de l'intégrale en résultent.  La dernière assertion sur la non-nullité résulte de  \cite[proposition  5.1]{FriJa}.
  \end{preuve}
\end{paragr}

\begin{paragr}[Période d'entrelacement $J_Q$.] ---  \label{S:JQ} Nous aurons besoin d'une variante de la construction ci-dessous. Rappelons qu'on a défini, au lemme \ref{lem:classes}, un élément $w_3\in \, _QW_H$. Pour alléger les notations, on pose $Q_3=Q_{w_3}$, $M_3=M_{Q_3}$ et $N_3=N_{Q_3}$ et  $Q_3^H=H\cap Q_3$, $M_3^H=H\cap M_3$ et $N_3^H=H\cap N_3$. On a la décomposition de Levi $Q_3^H=M_3^H N_3^H$. Les restrictions à $Z_{M_3}(\AAA)  M_3^H(F) N_{3}^H(\AAA)$ du caractère modulaire $\delta_{Q_3^H}$ du  groupe $ Q_3^H(\AAA)$ et du caractère
  \begin{align*}
    x\in Q_3(\AAA)\mapsto \exp(\bg \rho_Q^G,H_Q(w_3x)\bd)
  \end{align*}
coïncident.
  
  Pour tout $\varphi\in \Ac_{Q,\pi'}(G)$ on introduit  alors l'intégrale d'entrelacement
  \begin{align}\label{eq:periode entrelac2}
    J_Q(\varphi)=\int_{ [H]_{Q_3}  } \varphi(w_3 g)\, dg.
  \end{align}
  où l'on introduit $  [H]_{Q_3} =Z_{M_3}(\AAA)  M_3^H(F) N_{3}^H(\AAA)  \back H(\AAA) $.
  
\begin{proposition}\label{prop:entrelac2}
    \begin{enumerate}
    \item L'intégrale dans le membre de droite de \ref{eq:periode entrelac2} est absolument convergente et définit une forme linéaire $H(\AAA)$-invariante continue, notée $J_Q$, sur $\Ac_{Q,\pi'}(G)$.
    \item  La forme linéaire  $J_Q$ est non identiquement nulle sur $\Ac_{Q,\pi'}(G)$ si et seulement si la représentation $\sigma$ est symplectique.
    \end{enumerate}
  \end{proposition}

  \begin{preuve}  C'est une variante de la preuve de la proposition \ref{eq:periode entrelac}.  Soit  $g=(g',g_3)\in (G_{p+1}\times G_p)(\AAA)$. On utilise cette fois  la décomposition d'Iwasawa de $g'$ relativement au sous-groupe parabolique $P_{p,1}$ ce qui donne:
    \begin{align*}
      g'=  \begin{pmatrix}
        I_p & u \\ 0 & 1
      \end{pmatrix}
      \begin{pmatrix}
        g_1 & 0 \\ 0 &g_2
      \end{pmatrix} k
    \end{align*}
    avec $u\in \AAA^{p} $ identifié à une matrice colonne, $k\in K_{p+1}$, $g_2\in G_1(\AAA)$ et $g_1\in G_p(\AAA)$. On a alors
    \begin{align*}
      w_3 g= \begin{pmatrix}
        I_p & 0 &  u\\ 0 & I_p & 0 \\ 0 & 0 &1 
      \end{pmatrix}
      \begin{pmatrix}
        g_1 & 0  & 0\\ 0 &g_3 & 0 \\ 0 & 0 & g_2
      \end{pmatrix}w_3  \begin{pmatrix}
        k & 0 \\ 0 &I_p
      \end{pmatrix}.                                          
    \end{align*}
Il s'ensuit qu'on a 
    \begin{align*}
     \exp(\bg \rho_{Q}^{G}, H_{Q}(w_3g) \bd )    \delta_{Q_3^H}( \begin{pmatrix}
        g_1 & 0 & 0\\ 0 &g_2 &0 \\ 0 & 0 & g_3
      \end{pmatrix})^{-1}  =  |\det(g_3) \det(g_1)^{-1}|^{\frac12 }.
    \end{align*}
    On obtient
    \begin{align*}
      J_Q(\varphi)      &= \int_{ [ (G_p\times G_p)]_{G_{2p}} } \int_{K_{p+1}}\varphi_{-\rho_Q^G}\left(     \begin{pmatrix}
          g_1 & 0 &0\\ 0 &g_3& 0 \\ 0 & 0 & 1
      \end{pmatrix}  w_2 \begin{pmatrix}
          k & 0 \\ 0 &I_p
        \end{pmatrix} \right)\left|\frac{\det(g_3)}{\det(g_1)}\right|^{\frac12}dk dg_1dg_3.\\
    \end{align*}
    La fin de la preuve est alors similaire à celle de la proposition \ref{eq:periode entrelac}. 
  \end{preuve}
\end{paragr}

\begin{paragr}[Période d'entrelacement $J_{\bar P}$.] ---  \label{S:JPb} Soit $w_0$ l'unique élément du groupe de Weyl $W$  tel que $L=w_0 Mw_0^{-1}$ et $L\cap P_0= w_0 (M\cap P_0)w_0^{-1}$.  Notons qu'on a $w_0 \al=-\be$. Pour $\la\in \ago_{P,\CC}^{G,*}$, on dispose d'un opérateur d'entrelacement (pour les actions $I_P(\la)$ et $I_Q(w_0\la)$)
  \begin{align*}
    M(w_0, \la): \Ac_{P,\pi}\to \Ac_{Q,\pi'}
  \end{align*}
  qui est méromorphe en $\la$  et continu pour $\la$ dans le domaine d'holomorphie. Ce dernier comprend, du moins,   $i\ago_P^{G,*}$ et les $\la\in \ago_{P,\CC}^{G,*}$ de partie réelle  $\Re(\bg \la, \al^\vee\bd)$  assez grande, cf. \cite[corollaire 2.4]{Lap-remark}.  Il est unitaire (il préserve les normes de Petersson) pour $\la\in i\ago_P^{G,*}$.

  Soit $\bar P$ le sous-groupe parabolique de Levi $M$ tel que $P\cap \bar P=M$. Notons qu'on a $\bar P=w_{0}^{-1}Qw_0$. On définit, pour tout $\la\in \ago_{P,\CC}^{G,*}$, la   forme linéaire  $J_{\bar P} (\la)$ sur $\Ac_{P,\pi}(G)$ par
    \begin{align*}
    J_{\bar P} (\la,\varphi)=J_Q(M(w_0,\la)\varphi)
  \end{align*}
  pour tout $\varphi\in \Ac_{P,\pi}$. Notons que  $\la\mapsto J_{\bar P} (\la,\varphi)$ est une fonction méromorphe sur $\ago_{P,\CC}^{G,*}$ et holomorphe hors des hyperplans singuliers de $M(w_0,\la)$.
  \end{paragr}

\subsection{Périodes tronquées de séries d'Eisenstein}\label{ssec:period-Eis}

\begin{paragr}
Soit $\varphi\in \Ac_{P,\pi}(G)$, $\la\in  \ago_{P,\CC}^{G,*}$ et $g\in G(\AAA)$. On dispose alors de la série d'Eisenstein
  \begin{align*}
   E(g, \varphi, \la )= \sum_{\delta \in P(F)\back G(F)} \varphi_\la(\delta g).
  \end{align*}
 Cette égalité vaut lorsque la partie réelle de $\bg \la, \al^\vee\bd $ est assez grande:  la série dans le membre de droite est alors convergente et holomorphe en la variable $\la$. Le membre de gauche est ensuite défini par prolongement  méromorphe à $\ago_{P,\CC}^{G,*}$.  
\end{paragr}

\begin{paragr}   Nous aurons besoin du calcul des termes constants de $E( \varphi, \la )$ le long des sous-groupes paraboliques standard propres de $G$: ils  sont tous nuls hormis ceux le long des sous-groupes $P$ et $Q$ qui sont donnés par
 \begin{align}\nonumber
   E_P(g, \varphi,\la )&=\int_{[N_P]}E(ng, \varphi,\la )\,dn\\
 \label{eq:cst-P}      &=  \varphi_\la(g)
 \end{align}
 et
 \begin{align}\label{eq:cst-Q}
   E_Q(g, \varphi, \la )&= (M(w_0, \la)\varphi)_{w_0 \la}(g).
 \end{align}
 Rappelons que $w_0\in W$ a été défini au § \ref{S:JPb}.
 
\end{paragr}

\begin{paragr}[Série d'Eisenstein tronquée.] ---  Soit $\varphi\in \Ac_{P,\pi}(G)$ et $\la\in \ago_{P,\CC}^{G,*}$.  Soir $(\La^{T}_\theta  E)( \varphi, \la )$ la série d'Eisenstein $E(\varphi,\la)$ à laquelle on a appliqué l'opérateur de troncature $\La^{T}_\theta$ qui a été introduit au § \ref{S:intro-LaTtheta} et qui dépend d'un point auxiliaire $T\in \ago_{P_0}$. En tenant compte de sa définition  rappelée en \eqref{eq:intro-LaTtheta} et du calcul des termes constants effectué ci-dessus, on voit qu'on a 
  \begin{align}\nonumber
(\La^{T}_\theta  E)(g, \varphi, \la )   = E(g, \varphi,\la )-\sum_{w\in _P W_H  } \sum_{\delta\in (H\cap P_w)(F)\back H(F)} \hat\tau_P(H_P(w \delta g)-T_P)    \varphi_\la(wg)\\
\label{eq:LaTE}    -\sum_{w\in _Q W_H } \sum_{\delta\in (H\cap Q_w)(F)\back H(F)} \hat\tau_Q(H_Q(w \delta g)-T_Q)   (M(w_0, \la)\varphi)_{w_0 \la}(wg)
  \end{align}
  où $T_P$, resp.  $T_Q$, est le projeté de $T$ sur $\ago_P$, resp.  $\ago_Q$, parallèlement à $\Ker(\ago_{P_0}\to \ago_P)$, resp.  $\Ker(\ago_{P_0}\to \ago_Q)$. Notons que pour $H\in \ago_P$, resp. $\ago_Q$, on a $\hat\tau_P(H)=1$, resp.   $\hat\tau_Q(H)=1$, si et seulement si $\bg \al,H\bd >0$, resp.  $\bg \be,H\bd >0$. On pose
  \begin{align*}
      T_{\bar P}=w_0^{-1}T_Q.
  \end{align*}
  
  \begin{theoreme}\label{thm:LaTEis}
    Pour tout $\la\in \ago_{P,\CC}^{G,*}$ non singulier pour la série d'Eisenstein $E(\varphi, \la )$, l'intégrale
    \begin{align}\label{eq:LaTEis}
    \int_{[H]_G} (\La^{T}_\theta  E)(g, \varphi, \la ) \,dg.
    \end{align}
    est absolument convergente et est égale à
    \begin{align}
      \label{eq:calcul LaTEis}
      \frac{J_P(\varphi) \exp(\bg \la,T_P\bd)- J_{\bar P}(\la,\varphi) \exp(\bg \la,T_{\bar P}\bd)}{\theta_P(\la)}
    \end{align}
    où  $\theta_P(\la)$ est introduit en \eqref{eq:theta_P}.
  \end{theoreme}

  \begin{remarque}
    Il est usuel de définir la période régularisée de la série d'Eisenstein comme le terme constant de la partie polynomiale de la fonction de $T$ donnée par l'intégrale tronquée. Ici, on observe que,  pour des valeurs non nulles de $\la$, la période régularisée est nulle.
  \end{remarque}

  \begin{remarque}
    Friedberg-Jacquet ont défini une régularisation \emph{ad hoc} de la période de la série d'Eisenstein qu'on considère, cf. \cite[théorème 5.1]{FriJa}, qui correspond ici à la période d'entrelacement $J_P$. Le théorème \ref{thm:LaTEis} donne donc une autre façon de faire apparaître la distribution $J_P$.
  \end{remarque}
  
  La démonstration du théorème \ref{thm:LaTEis} se trouve au § \ref{S:preuve LaTEis}. Auparavant, nous donnons deux corollaires au théorème.
\end{paragr}

\begin{paragr}[Un premier corollaire] --- Celui-ci énonce  une équation fonctionnelle pour les périodes d'entrelacement.

  \begin{corollaire}\label{cor:LaTEis}
    Pour tout  $\varphi\in \Ac_{P,\pi}(G)$, on  a
    \begin{align}\label{eq:eq fct}
      J_{\bar P}(0,\varphi) =J_P(\varphi)
    \end{align}
    ou encore, par définition du membre de gauche,
    \begin{align*}
      J_Q(M(w_0,0)\varphi)=J_P(\varphi).
    \end{align*}
  \end{corollaire}

  \begin{remarque}
  La condition   \ref{eq:eq fct} implique que le couple de fonctions $(J_P(\varphi),  J_{\bar P}(\la,\varphi))$ de la variable $\la\in \ago_{P,\CC}^{G,*}$ (la première fonction est en fait constante) est une $(G,M)$-famille au sens d'Arthur, cf. \cite{arthur2}.
  \end{remarque}
  
  \begin{preuve}
    La série d'Eisenstein $E(\varphi, \la )$ est holomorphe pour $\la$ dans un voisinage de $0$. Il s'ensuit que  l'intégrale \eqref{eq:LaTEis} est également holomorphe dans un voisinage de $0$. Le théorème \ref{thm:LaTEis} implique que l'expression \eqref{eq:calcul LaTEis} est holomorphe en $0$. En particulier, son dénominateur doit s'annuler en $0$ d'où \eqref{eq:eq fct}.
  \end{preuve}
\end{paragr}

\begin{paragr}[Dérivée d'opérateurs d'entrelacements.] ---  \label{S:M derive} On définit  un opérateur continu
  \begin{align}\label{eq:mc}
    \mc_{P,\pi}:\Ac_{P,\pi}(G)\to \Ac_{P,\pi}(G)
  \end{align}
   par la formule
 \begin{align}\label{eq:formule-mc}
       \mc_{P,\pi}  \varphi=\lim_{\la\to 0}   \frac{\varphi-M(w_0,0)^{-1}  M(w_0,\la)\varphi}{\theta_P(\la)}
  \end{align}
  pour tout $\varphi\in \Ac_{P,\pi}(G)$. On renvoie à  \eqref{eq:theta_P} pour la définition de $\theta_P(\la)$. L'existence de la limite et la continuité sont une conséquence immédiate de l'holomorphie et de l'uniforme continuité de l'opérateur   $M(w_0,\la)$ au voisinage de $\la=0$. Cet opérateur n'est qu'un exemple, parmi d'autres, d'opérateurs introduits par Arthur via la théorie des $(G,M)$-familles qui interviennent dans la définition des caractères pondérés d'Arthur,   cf. \cite[section 15, éq. (15.14)]{Ar-cours}. C'est un opérateur qui n'est pas équivariant pour l'action par translation à gauche de $G(\AAA)$. Pour formuler ce défaut d'équivariance,  on introduit, pour tout $y\in G(\AAA)$ un opérateur $I_P'(y): \Ac_{P,\pi}(G)\to \Ac_{P,\pi}(G)$ défini pour tout $\varphi\in \Ac_{P,\pi}(G)$ et tout $x\in G(\AAA)$
  \begin{align}\label{eq:IP'}
    (I_P'(y)\varphi)(x)=  \lim_{\la\to 0}   \frac{    (I_P(\la,y)\varphi-I_P(y)\varphi)(x)  }{\theta_P(\la)},
  \end{align}
la limite ci-dessus étant prise sur $\la\in \ago_{P,\CC}^{G,*}$.  Introduisons la fonction
\begin{align}
  \label{eq:fct-kP}
k_P:P(\AAA)\back G(\AAA)\to (K\cap P(\AAA))\back K
\end{align}
définie par la condition $x\in P(\AAA) k_P(x)$ pour tout $x\in G(\AAA)$. Un calcul élémentaire montre que, pour tout $\varphi\in \Ac_{P,\pi}(G)$, on a
\begin{align}\label{eq:IP'bis}
   (I_P'(y)\varphi)(x)=    \frac{  \bg \al, H_P( k_P(x) y)\bd}{\theta_P(\al)}\varphi(xy).  
\end{align}
Lorsqu'on remplace $P$ par $Q$ et $\al$ par $\be$, on obtient un opérateur $I_Q'(y): \Ac_{Q,\pi'}(G)\to \Ac_{Q,\pi'}(G)$. Les opérateurs $\mc_{P,\pi}$ et $ I_P(y)$ ne commutent pas en général comme le montre le lemme suivant.

\begin{lemme}\label{lem:crochet}
    Pour tout $y \in G(\AAA)$, on a l'égalité suivante
    \begin{align*}
      \mc_{P,\pi} I_P(y) - I_P(y) \mc_{P,\pi}= I_P'(y) +   M(w_0,0)^{-1} I_Q'(y) M(w_0,0).
    \end{align*}
  \end{lemme}

  \begin{preuve}
    Soit $\varphi\in \Ac_{P,\pi}(G)$ et $y \in G(\AAA)$. En utilisant la définition \eqref{eq:formule-mc} de l'opérateur $\mc_{P,\pi}$, on obtient
    \begin{align*}
      M(w_0,0)  (\mc_{P,\pi} I_P(y) \varphi - I_P(y) \mc_{P,\pi}\varphi )=\lim_{\la\to 0}   \frac{I_Q(y) M(w_0,\la)\varphi-  M(w_0,\la)I_P(y)\varphi}{\theta_P(\la)}.
    \end{align*}
    Cette limite est la somme des limites des termes 
    \begin{align*}
      \frac{I_Q(y) -  I_Q(w_0\la,y) }{\theta_P(\la)}M(w_0,\la)\varphi
    \end{align*}
    et
    \begin{align*}
    M(w_0,0)  \frac{I_P(\la,y)-  I_P(y)}{\theta_P(\la)} \varphi.
    \end{align*}
    Le résultat est alors évident.
  \end{preuve}
\end{paragr}

\begin{paragr}[Un second corollaire.] --- Celui-ci explicite le théorème  \ref{thm:LaTEis} au point $\la=0$ et fait intervenir  l'opérateur $\mc_{P,\pi}$ défini en \eqref{eq:mc}.

\begin{corollaire} \label{cor:LaTEisven 0}Pour tout $\varphi\in \Ac_{P,\pi}(G)$, on a  
  \begin{align}\label{eq:LaTEis en 0}
    \int_{[H]_G} (\La^{T}_\theta  E)(g, \varphi, 0) \,dg= J_P(    \mc_{P,\pi}   \varphi) +  J_{P}(\varphi)   \bg \al, T_P-T_{\bar P}\bd \theta_P(\al)^{-1} .
  \end{align}
\end{corollaire}

\begin{remarque}
  Si l'on définit la période régularisée de la série d'Eisenstein $E(g, \varphi, 0)$ comme le terme constant de la fonction de $T$ donnée par l'intégrale tronquée, on observe que  la période régularisée  ainsi définie est donnée par la distribution $\varphi\mapsto J_P(    \mc_{P,\pi}   \varphi) $. Celle-ci  n'est pas invariante comme le laisse supposer la formule du lemme \ref{lem:crochet}: cela reflète le fait que la partie polynomiale en $T$ de  l'intégrale tronquée (ici l'intégrale tronquée coïncide avec un polynôme de degré $1$) n'est pas constante.
\end{remarque}

\begin{preuve}
  Comme la série d'Eisenstein est holomorphe en $\la=0$, il en est de même de l'expression \eqref{eq:calcul LaTEis}. Le théorème \ref{thm:LaTEis} implique que le membre de gauche de \eqref{eq:LaTEis en 0} est donné par la limite de l'expression
  \begin{align*}
      \frac{J_P(\varphi) \exp(\bg \la,T_P\bd)- J_{\bar P}(\la,\varphi) \exp(\bg \la,T_{\bar P}\bd)}{\theta_P(\la)},
  \end{align*}
  quand $\la \in \ago_{P,\CC}^{G,*}$ tend vers $0$  et  $\bg \la,\al^\vee\bd\not=0$.  L'expression ci-dessus est la somme de

     \begin{align}\label{eq:lim1}
      \frac{J_P(\varphi) - J_{\bar P}(\la,\varphi) }{\theta_P(\la) }\exp(\bg \la,T_{P}\bd)
     \end{align}
    et \begin{align}\label{eq:lim2}
      J_{\bar P}(\la,\varphi)   \frac{ \exp(\bg \la,T_P\bd)- \exp(\bg \la,T_{\bar P}\bd)}{\theta_P(\la)}.
    \end{align}
    L'opérateur $M(w_0,\la)$ est holomorphe en $\la=0$ de même que  $J_{\bar P}(\la,\varphi)$. 
    En utilisant l'équation fonctionnelle du corollaire \ref{cor:LaTEis}, on obtient que la limite de \eqref{eq:lim2}, quand $\la$ tend vers $0$, est égale à
    \begin{align*}
        J_{P}(\varphi)   \bg \al, T_P-T_{\bar P}\bd \theta_P(\al)^{-1}.
    \end{align*}
    Pour obtenir la limite de \eqref{eq:lim1}, on écrit, toujours en utilisant  l'équation fonctionnelle du corollaire \ref{cor:LaTEis}:
    \begin{align*}
      J_P(\varphi) - J_{\bar P}(\la,\varphi) &=J_P(\varphi)-J_Q(M(w_0,\la)\varphi)\\
      &= J_P(\varphi-M(w_0,0)^{-1}M(w_0,\la)\varphi).
    \end{align*}
    En utilisant la continuité de la forme linéaire $J_P$ et la définition \eqref{eq:mc} de l'opérateur   $\mc_{P,\pi}$,  on a
    \begin{align*}
      \frac{J_P(\varphi) - J_{\bar P}(\la,\varphi) }{\theta_P(\la) }&=J_P(  \frac{\varphi-M(w_0,0)^{-1}  M(w_0,\la)\varphi}{\theta_P(\la)})\\
      &\to_{\la \to 0}  J_P(  \mc_{P,\pi}\varphi).
    \end{align*}
    La conclusion est claire.
\end{preuve}
\end{paragr}

\begin{paragr}[Démonstration du théorème \ref{thm:LaTEis}.] ---  \label{S:preuve LaTEis} On reprend et on étend les notations du lemme \ref{lem:classes} et de la section \ref{ssec:entrelacement}. Pour $i=1,2$, on note $P_i=P_{w_i}$ et $P_i^H=P_{w_i}\cap H$. On commence par observer en utilisant le  lemme \ref{lem:classes} et l'expression \eqref{eq:LaTE} que l'intégrale \eqref{eq:LaTEis} est (formellement) égale à la somme des cinq intégrales:
    \begin{align}
\label{eq:I0}      &  I_0^T(\varphi,\la)=\int_{Z_G(\AAA) P^H(F) \back H(\AAA)  } \left(1- \hat\tau_P(H_P(g)-T_P) \right) \varphi_\la ( g)\, dg\,;\\
  \label{eq:I1}     &   I_1^T(\varphi,\la)=\int_{Z_G(\AAA) P_{1}^H(F) \back H(\AAA)} \left(1- \hat\tau_P(H_P(w_1      g)-T_P)  \right)\varphi_\la (w_1  g)\, dg\,;\\
      \label{eq:I2}     & I_2(\varphi,\la)= \int_{Z_G(\AAA) P_{2}^H(F)\back H(\AAA)} \varphi_\la (w_2 g)\, dg \,;\\
      \label{eq:I3}     & I_3^T(\varphi,\la)= -\int_{Z_G(\AAA)  Q^H_{3}(F) \back H(\AAA)}  \hat\tau_Q(H_Q(w_3g)-T_Q)  (M(w_0,\la)\varphi)_{w_0\la} ( w_3g)\, dg \,;\\
    \label{eq:I4}     & I_4^T(\varphi,\la)= - \int_{Z_G(\AAA)  Q^H(F) \back H(\AAA)}   \hat\tau_Q(H_Q( g)-T_Q)  (M(w_0,\la)\varphi)_{w_0\la} ( g)\, dg.
    \end{align}

    La proposition suivante implique la validité du théorème \ref{thm:LaTEis} sur un ouvert. Par prolongement méromorphe, on obtient l'énoncé général.
    
  \begin{proposition}\label{prop:preuve LaTEis}
    Pour $0\leq k\leq 4$ et tout $\la\in \ago_{P,\CC}^{G,*}$ tel que $\Re(\bg\la,\al^\vee\bd)$ assez positif,  les intégrales ci-dessus sont absolument convergentes. De plus, elle sont toutes nulles sauf éventuellement dans les deux cas suivants pour lesquels on a:
    \begin{align*}
      & I_0^T(\varphi,\la)= J_P(\varphi)       \exp(\bg \la,T\bd) \theta_P(\la)^{-1}\\
      &I_3^T(\varphi,\la)=-J_{\bar P}( M(w_0,\la)\varphi) \exp(\bg\la,T_{\bar P}\bd) \theta_P(\la)^{-1}.
    \end{align*}
  \end{proposition}

  Le reste de la sous-section est consacrée à la démonstration de la proposition \ref{prop:preuve LaTEis} qui s'inspire d'ailleurs très fortement de \cite[section 5]{FriJa}.
\end{paragr}

\begin{paragr}[Calcul de $ I_0^T(\varphi,\la)$ et $I_3^T(\varphi,\la)$.] ---   Rappelons qu'on a  $P^H=M^H N_P^H$. En tenant compte des observations du § \ref{S:JP} sur les caractères modulaires, on voit que  $I_0(\varphi,\la)$ est égal à 
  \begin{align*}
    &\int_{  Z_M(\AAA) M^H(F) N_{P}^H(\AAA)  \back H(\AAA)  }  \int_{Z_G(\AAA) \back Z_M(\AAA)}  \left(1- \hat\tau_P(H_P(a g)-T_P) \right)  \exp(\bg -\rho_P^G,H_P(a)\bd    \varphi_\la(a g)\, da dg\\
    & \int_{  Z_M(\AAA) M^H(F) N_{P}^H(\AAA)   \back H(\AAA)  }  \varphi_\la(g)  \int_{\ago_P^G}  \left(1- \hat\tau_P(X+H_P(g)-T_P) \right)  \exp(\bg \la,X\bd)  \, dX \,dg.
  \end{align*}
 L'intégrale intérieure est convergente pour $\Re(\bg \la, \al^\vee\bd)>0$ et vaut:
  \begin{align*}
    \exp(\bg \la,T_P- H_P(g)\bd) \theta_P(\la)^{-1}.
  \end{align*}
  Le résultat pour $ I_0^T(\varphi,\la)$ s'ensuit compte tenu de la proposition \ref{prop:periode entrelac}.

  Traitons de la même façon $I_3^T(\varphi,\la)$. En tenant compte des  caractères modulaires, cf. § \ref{S:JQ}, on obtient que $-I_3(\varphi,\la)$ est égal à l'intégrale sur $g\in  Z_{M_3}(\AAA) M^H_3(F)  N_{3}^H(\AAA)  \back H(\AAA) $ de 
  \begin{align*}
    & \int_{Z_G(\AAA) \back Z_{M_3}(\AAA)} \hat\tau_Q(H_Q(w_3 a g)-T_Q) \exp(\bg -\rho_Q^G,H_Q(w_3aw_3^{-1})\bd )  (M(w_0,\la)\varphi)_{w_0\la} ( w_3a g) \, da \\
    &=  (M(w_0,\la)\varphi)_{w_0\la} ( w_3g) \int_{\ago_Q^{G}} \hat\tau_Q(X+H_Q(w_3 g)-T)_Q \exp(\bg w_0\la,X \bd  ) \, dX .
  \end{align*}

  L'intégrale ci-dessous converge pour $\Re(\bg w_0\la,\be^\vee\bd)  <0$ c'est-à-dire $\Re(\bg \la,\al^\vee\bd)  >0$ et vaut
   \begin{align*}
    -\exp(\bg w_0\la,T_Q- H_Q(w_3g)\bd) \vol(\ago_Q^G / \ZZ \be^\vee) \bg w_0\la, \be^\vee\bd^{-1}.
  \end{align*}
On a   $\theta_P(\la)^{-1}=-\vol(\ago_Q^G / \ZZ \be^\vee) \bg w_0\la, \be^\vee\bd^{-1}$. On conclut alors  comme auparavant en utilisant cette fois la  proposition \ref{prop:entrelac2}.
\end{paragr}

\begin{paragr}[ Calcul de $ I_1^T(\varphi,\la)$.] --- \label{S:calcul I1}Commençons par des manipulations formelles que l'on justifie ensuite. Soit $M_1=w_1^{-1}Mw_1$ et $N_1=N_{P_{1}}$. On note avec un exposant $H$ leur intersection avec $H$. Soit $\delta_{P_1^H}$ le caractère modulaire de $P_1^H(\AAA)$. L'intégrale $ I_1^T(\varphi,\la)$ s'écrit comme l'intégrale sur $g\in  Z_{M_1}(\AAA) M_1^H(F) N_{1}^H(\AAA)  \back H(\AAA) $ de 
  \begin{align}\label{eq:I_1 step1}
 \int_{Z_G(\AAA) \back Z_{M_1}(\AAA)} \left(1-\hat\tau_P(H_P(w_1 a g)-T_P) \right)  \delta_{P_1^H}(a)^{-1}   \varphi_{\la} ( w_1a g) \, da .
  \end{align}

  On écrit la décomposition d'Iwasawa  de $g\in H(\AAA)$ relative au sous-groupe parabolique $P_1^H$ sous la forme suivante:
  \begin{align*}
    g=  \begin{pmatrix}
      I_{p+1} & 0 & 0 \\ 0 & 1 & u \\ 0 & 0 & I_{p-1}
      \end{pmatrix} \begin{pmatrix}
      g_1 & 0 & 0 \\ 0 & g_2 & 0 \\ 0 & 0 & g_3
      \end{pmatrix} \begin{pmatrix}
      I_{p+1} & 0  \\ 0 & k 
    \end{pmatrix}
  \end{align*}
     avec $u\in \AAA^{p-1} $ identifié à une matrice ligne, $k\in K_{p}$, $g_1\in G_{p+1}(\AAA)$, $g_2\in G_1(\AAA)$ et $g_3\in G_{p-1}(\AAA)$. 
             On a alors
             \begin{align}\label{eq:I1 rho}
   \exp(\bg \rho_P^G,H_P(w_1g)\bd)        \delta_{P_1^H}( \begin{pmatrix}
      g_1 & 0 & 0 \\ 0 & g_2 & 0 \\ 0 & 0 & g_3
      \end{pmatrix}   )^{-1} =|g_2|\cdot   \left|\frac{\det(g_3)}{\det(g_1)}\right|^{1/2}.
             \end{align}
             On en déduit aussi que,  pour tout $a\in Z_M(\AAA)$, on a
             \begin{align*}
              \exp( \bg \rho_P^G,H_P(a)\bd)        \delta_{P_1^H}(w_1^{-1}aw_1)^{-1} = \exp(\bg \al, H_P(a)\bd).
             \end{align*}
             On en déduit, par un changement de variables, que l'intégrale \eqref{eq:I_1 step1} est égale à
\begin{align*}
  \varphi_{\la} ( w_1g)\int_{\ago_P^G} \left(1-\hat\tau_P(X+H_P(w_1 g)-T_P) \right)\exp(\bg \al+\la, X \bd)   \, dX. 
  \end{align*}
  L'intégrale ci-dessus converge absolument pour $\Re(\bg \al+\la, \al^\vee\bd)>0$ et elle  vaut sur ce domaine
  \begin{align*}
    \exp(\bg \al+\la, T_P-  H_P(w_1 g)\bd) \theta_P(\al+\la)^{-1}.
  \end{align*}
  On est alors ramené à prouver la convergence absolue et obtenir le calcul de  l'intégrale suivante:
  \begin{align*}
    \int_{  Z_{M_1}(\AAA) M_1^H(F) N_{1}^H(\AAA)  \back H(\AAA) }\varphi(w_1g) \exp(-\bg \al, H_P(w_1 g)\bd) \, dg.
  \end{align*}
  
  Vu \eqref{eq:I1 rho}, par décomposition d'Iwasawa, on est ramené à l'intégrale suivante:
  
  \begin{align*}
         \int_{ [G_{p+1}\times G_{p-1}]_{G_{2p}}} \int_{K_{p}}\varphi_{-\rho_P^G}\left(     \begin{pmatrix}
          1 & 0 &0\\ 0 &g_1& 0 \\ 0 & 0 & g_3
      \end{pmatrix}  w_1  \begin{pmatrix}
          I_{p+1} & 0 \\ 0 &k
      \end{pmatrix} \right)  \frac{|\det(g_3)|^{\frac12+\frac1{2p}}}{|\det(g_1)|^{\frac12-\frac1{2p}}} dk dg_1dg_3.
  \end{align*}
Comme $\pi$ est cuspidale, l'intégrande est à décroissance rapide en $(g_1,g_3)\in [G_{p+1}\times G_{p-1}]_{G_{2p}}$ d'où la  convergence absolue.  L'intégrale est en fait nulle. Il suffit de voir que pour toute représentation automorphe cuspidale $\sigma$ de $G_{2p}(\AAA)$ et toute fonction $\varphi\in \Ac_\sigma(G_{2p})$, on a
  \begin{align*}
      \int_{ [G_{p+1}\times G_{p-1}]_{G_{2p}}}  \varphi\left(     \begin{pmatrix}
        g_1& 0 \\ 0 & g_3
      \end{pmatrix}  \right)  \frac{|\det(g_3)|^{\frac12+\frac1{2p}}}{|\det(g_1)|^{\frac12-\frac1{2p}}} dg_1dg_3=0.
  \end{align*}
  Si $p=1$, l'intégrale n'est autre que le produit scalaire d'une fonction automorphe cuspidale sur $[G_2]$ contre une fonction constante: elle est donc nulle. Si $p>1$, l'intégrale est aussi nulle en vertu de la \cite[proposition 2.1]{FriJa}.
\end{paragr}

\begin{paragr}[ Calcul de $ I_4^T(\varphi,\la)$.] --- Ici encore, nous commençons par des manipulations formelles qui seront justifiées \emph{in fine}. Rappelons que $L$ est le facteur de Levi standard de $Q$. On a une décomposition de Levi $Q^H=L^H N_Q^H$ avec $L^H=L\cap H$ et $N_Q^H=N_Q\cap H$. Soit $\delta_{Q^H}$ le caractère modulaire du groupe $Q^H(\AAA)$. On écrit l'intégrale $ -I_4^T(\varphi,\la)$  comme l'intégrale sur $g\in  Z_L(\AAA) L^H(F) N_{Q}^H(\AAA)  \back H(\AAA) $ de 
  \begin{align}\label{eq:I_4 step1}
 \int_{Z_G(\AAA) \back Z_L(\AAA)}\hat\tau_Q(H_Q(a g)-T_Q) \delta_{Q^H}(a)^{-1} (M(w_0,\la)\varphi)_{w_0\la} ( a g) \, da .
  \end{align}
  On écrit la décomposition d'Iwasawa  de $g\in H(\AAA)$ relative au sous-groupe parabolique $Q^H$:
  \begin{align*}
    g=  \begin{pmatrix}
      I_{p+1} & 0 & 0 \\ 0 &  I_{p-1}& u \\ 0 & 0 & 1
      \end{pmatrix} \begin{pmatrix}
      g_1 & 0 & 0 \\ 0 & g_2 & 0 \\ 0 & 0 & g_3
      \end{pmatrix} \begin{pmatrix}
      I_{p+1} & 0  \\ 0 & k 
    \end{pmatrix}
  \end{align*}
     avec $u\in \AAA^{p-1} $ identifié à une matrice colonne, $k\in K_{p}$, $g_1\in G_{p+1}(\AAA)$, $g_2\in G_{p-1}(\AAA)$ et $g_3\in G_1(\AAA)$.   On calcule:
             \begin{align}\label{eq:I4 rho}
\exp( \bg \rho_Q^G,H_Q(g)\bd)  \delta_{Q^H}(  \begin{pmatrix}
      g_1 & 0 & 0 \\ 0 & g_2 & 0 \\ 0 & 0 & g_3
      \end{pmatrix} )^{-1}  = |g_3|^{-1}\cdot   \left|  \frac{ \det(g_1)}{\det(g_2)}\right|^{1/2} .
             \end{align}
             En particulier, pour tout $a\in Z_L(\AAA)$, on a
             \begin{align}
                    \label{eq:I4 rho2}          \exp( \bg \rho_Q^G,H_Q(a)\bd   )   \delta_{Q^H}(  a)^{-1}=\exp( \bg \be,H_Q(a)\bd).
             \end{align}
             Comme au § \ref{S:calcul I1}, on voit que l'intégrale \eqref{eq:I_4 step1} est égale à l'intégrale ci-dessous qui converge pour $\Re(\bg \be +w_0\la,\be^\vee\bd)<0$ c'est-à-dire $\Re(\bg \la -\al,\al^\vee\bd)>0$ 
             \begin{align*}
               (M(w_0,\la)\varphi)_{w_0\la} ( g)       \int_{\ago_L^G}\hat\tau_Q(X+H_Q( g)-T_Q) \exp(\bg \be +w_0\la,X\bd) \, dX\\
               = (M(w_0,\la)\varphi)_{w_0\la} ( g)    \exp(\bg \be +w_0\la,T_Q-H_{Q}(g)\bd)  \theta_P( \la-\al)^{-1}.
             \end{align*}
             Pour conclure, on observe que
             \begin{align*}
               \int_{ Z_L(\AAA) L^H(F) N_{Q}^H(\AAA)   \back H(\AAA)} (M(w_0,\la)\varphi)_{w_0\la} ( g)    \exp(-\bg \be +w_0\la,H_{Q}(g)\bd)\, dg\\
               =\int_{ [G_{p+1}\times G_{p-1}]_{G_{2p}}} \int_{K_{p}}\varphi_{-\rho_Q^G}\left(     \begin{pmatrix}
          g_1& 0 &0\\ 0 &g_2& 0 \\ 0 & 0 & 1
      \end{pmatrix}  w_1  \begin{pmatrix}
          I_{p+1} & 0 \\ 0 &k
      \end{pmatrix} \right)  \frac{|\det(g_1)|^{\frac12+\frac1{2p}}}{|\det(g_2)|^{\frac12-\frac1{2p}}} dk dg_2dg_1.
             \end{align*}
             Pour les mêmes raisons qu'au § \ref{S:calcul I1}, cette intégrale est absolument convergente et nulle.
\end{paragr}

  \begin{paragr}[Calcul de $ I_2(\varphi,\la)$.] ---  On  a
    \begin{align*}
       w_2^{-1}=
      \begin{pmatrix}
        - 1 & 1 & 0 \\ 0 & 0 &  I_{2p-1}\\ 1 & 0&0
      \end{pmatrix}.
    \end{align*}
et $P_{2}^H$ est le groupe des matrices de la forme
    \begin{align*}
 \begin{pmatrix}
        a  & u & 0  & 0 \\ 0 & g_2 &  0 & 0 \\ 0 &0 & g_3& 0 \\  0 &0 & v & a 
      \end{pmatrix}
    \end{align*}
    avec $a\in G_1$, $g_2\in G_p$, $g_3\in G_{p-1}$, $u$ et $v$ des vecteurs ligne de taille respective $p$ et $p-1$. Soit $g=(g',g'')\in (G_{p+1}\times G_p)(\AAA)$. On  va utiliser  la décomposition d'Iwasawa de $g$ relativement au sous-groupe parabolique $P_{1,p}\times \bar P_{p-1,1}$ : pour cela, on écrit
    \begin{align*}
     & g'=       \begin{pmatrix}
        1  & u  \\ 0  & I_p
      \end{pmatrix}     \begin{pmatrix}
        g_1  & 0  \\ 0  & g_2      \end{pmatrix} k_1\\
     &  g''=       \begin{pmatrix}
        I_{p-1}  & 0  \\ v  & 1
      \end{pmatrix}\begin{pmatrix}
        g_3  & 0 \\ 0  & g_4
      \end{pmatrix} k_2
          \end{align*}
    avec $k_1\in K_{p+1}$, $k_2\in K_p$, $g_1,g_4\in G_1(\AAA)$, $g_2\in G_p(\AAA)$, $g_3\in G_{p-1}(\AAA)$, $u$ et $v$ des vecteurs ligne de taille respective $p$ et $p-1$. 

    On a alors
    \begin{align*}
&      w_2  \begin{pmatrix}
        1 & u & 0  & 0 \\ 0 & I_{p} &  0 & 0 \\ 0 &0 & I_{p-1} & 0 \\  0 &0 & v & 1 
      \end{pmatrix}w_2^{-1}= \begin{pmatrix}
        1 & 0& 0  & v \\ 0& 1 &  u & v \\ 0 &0 & I_p & 0 \\  0 &0 & 0 & I_{p-1} 
      \end{pmatrix}\in N_P(\AAA)    \begin{pmatrix}
        1 & 0& 0  & 0 \\ 0& 1 &  u & v \\ 0 &0 & I_p & 0 \\  0 &0 & 0 & I_{p-1} 
      \end{pmatrix}\\
 &      w_2  \begin{pmatrix}
        g_1 & 0 & 0  & 0 \\ 0 & g_2 &  0 & 0 \\ 0 &0 & g_3 & 0 \\  0 &0 & 0 & g_4
      \end{pmatrix}w_2^{-1}= \begin{pmatrix}
        g_4& 0& 0  & 0 \\ g_4-g_1 & g_1 &  0 & 0 \\ 0 &0 & g_2 & 0 \\  0 &0 & 0 & g_3 
      \end{pmatrix}.
    \end{align*}
    À l'aide de  la décomposition d'Iwasawa, on voit que l'intégrale $I_2(\la,\varphi)$ est formellement égale à
        \begin{align}\label{eq:calcul-I2}
\int |g_1 ^{-p}g_4^{-(p-1)}|\cdot   |\det(g_2)\det(g_3)|  \varphi_\la\left( \begin{pmatrix} 1 & 0 \\ 0 & n  
      \end{pmatrix} \begin{pmatrix}
        g_4 & 0& 0  & 0 \\ g_4-g_1 & g_1 &  0 & 0 \\ 0 &0 & g_2 & 0 \\  0 &0 & 0 & g_3
      \end{pmatrix} w_2k      \right)   \, dn d\underline{g}  dk
    \end{align}
    où l'intégrale est prise sur $n\in  [N_{1,2p-1}] $, $\underline{g} =(g_1,g_2,g_3,g_4)\in [G_1\times G_p\times G_{p-1}\times G_1   ]_G$ et $k\in K_{p+1}\times K_p$ (ce dernier groupe est naturellement vu comme sous-groupe de $K_{2p+1}$).
    On effectue le changement de variables $(g_1,g_4)\mapsto (g_1,g_1^{-1}g_4)$ ce qui donne 
      \begin{align*}
\int |g_1 ^{-2p+1}g_4^{-p+1}|\cdot   |\det(g_2)\det(g_3)| \varphi_\la\left( \begin{pmatrix} 1 & 0 \\ 0 & n  
      \end{pmatrix} \begin{pmatrix}
        g_4g_1 & 0& 0  & 0 \\ g_4g_1-g_1 & g_1 &  0 & 0 \\ 0 &0 & g_2 & 0 \\  0 &0 & 0 & g_3
      \end{pmatrix} w_2k      \right)   \, dn d\underline{g}' dg_4 dk
      \end{align*}
      où l'on intègre sur $n$ et $k$ sont comme ci-dessus, $g_4\in [G_1]$ et $\underline{g}' =(g_1,g_2,g_3)\in [G_1\times G_p\times G_{p-1}  ]_{G_{2p}}$.
      On veut établir la convergence absolue de cette intégrale.     On peut et on va supposer que $n$ reste dans un compact fixé de $N_{1,2p-1}(\AAA)$.  On peut aussi écrire $g_4= x t $ avec $t\in \RR_+^\times$ identifié à un idèle par le choix d'une place archimédienne $w$ de $F$ et $x$ qui reste dans une partie compacte fixée du sous-groupe des idèles de norme $1$. On a alors
        \begin{align*}
        \begin{pmatrix}
          g_4 g_1& 0 \\ g_4g_1-g_1& g_1
        \end{pmatrix}&=  \begin{pmatrix}
          t & 0 \\  t& 1 
        \end{pmatrix}   \begin{pmatrix}
          x & 0 \\  0& 1 
        \end{pmatrix} \begin{pmatrix}
          g_1 & 0 \\ 0& g_1
        \end{pmatrix}   \begin{pmatrix}
          1 & 0 \\ -1& 1
        \end{pmatrix}\\
&=           \begin{pmatrix}
          t & 0 \\  t& 1 
        \end{pmatrix}  \begin{pmatrix}
          g_1 & 0 \\ 0& g_1
        \end{pmatrix}   \begin{pmatrix}
          x & 0 \\ -1& 1
        \end{pmatrix}.
        \end{align*}
       La décomposition d'Iwasawa s'écrit
\begin{align*}
       \begin{pmatrix}
          t & 0 \\  t& 1 
       \end{pmatrix}=\begin{pmatrix}
          1 & v_{t} \\  0& 1 
        \end{pmatrix} \begin{pmatrix}
          t y_t^{-1} & 0 \\  0& y_t
        \end{pmatrix}k_t
\end{align*}

      où $v_{t}=t^{2}(t^2+1)^{-1}$, $y_t=\sqrt{t^2+1}$ et $k_t\in K_{2,w}\subset G_2(F_w)$. On a donc aussi
      \begin{align*}
     \begin{pmatrix}
        g_4g_1 & 0& 0  & 0 \\ g_4g_1-g_1 & g_1 &  0 & 0 \\ 0 &0 & g_2 & 0 \\  0 &0 & 0 & g_3
      \end{pmatrix}   =  
\begin{pmatrix}
          1 & v_{t}   &0 & 0 \\  0& 1  &0 & 0\\0 &0 & I_p & 0 \\  0 &0 & 0 & I_{p-1}
        \end{pmatrix}
        \begin{pmatrix}
        g_1 & 0& 0  & 0 \\ 0 & g_1 &  0 & 0 \\ 0 &0 & g_2 & 0 \\  0 &0 & 0 & g_3
      \end{pmatrix} \times \\\begin{pmatrix}
          t y_t^{-1} & 0& 0  & 0 \\ 0 & y_t &  0 & 0 \\ 0 &0 & I_p & 0 \\  0 &0 & 0 & I_{p-1}
      \end{pmatrix} k_t\begin{pmatrix}
        x& 0& 0  & 0 \\ -1 & 1 &  0 & 0 \\ 0 &0 & I_p & 0 \\  0 &0 & 0 & I_{p-1}
      \end{pmatrix} 
      \end{align*}
      On identifie $\la\in \ago_{P,\CC}^{G,*}$ à $2s\rho_P^G$   pour  un certain $s\in \CC$. Pour vérifier la convergence, on va prendre $s\in \RR$.

      Finalement, en utilisant le changement de variables $g_1\mapsto g_1 y_t$ et la définition de $y_t$, on est  ramené à prouver la convergence   de l'intégrale 
      \begin{align*}
        \int  (|g_1|^{2p-1}  |\det(g_2)\det(g_3)|^{-1})^{s-1/2}     \frac{t^{2ps+1}}{(  t^2+1)^{(p+1/2)(s+1/2)}} \\
        \left|\varphi_{-\rho}\left(   \begin{pmatrix} 1 & 0 \\ 0 & n  
      \end{pmatrix}   \begin{pmatrix}
1 & 0& 0  & 0 \\ 0 & g_1 &  0 & 0 \\ 0 &0 & g_2 & 0 \\  0 &0 & 0 & g_3
      \end{pmatrix}  k'\right) \right|\,dt dg_1 dg_2 dg_3,
      \end{align*}
      où l'on intègre sur $t\in   \RR_+^\times$ (pour une mesure de Haar $dt$ sur    $\RR_+^\times$) et   $(g_1,g_2,g_3)\in [G_1\times G_p\times G_{p-1}  ]_{G_{2p}}$,   et à la majorer uniformément pour $n$ et $k'$ restant dans des compacts fixés.  

      Il résulte de la cuspidalité de $\sigma$ que, pour tout $N>0$, il  existe $C>0$ de sorte que 
      \begin{align*}
         \left|\varphi_{-\rho}\left(   \begin{pmatrix} 1 & 0 \\ 0 & n  
      \end{pmatrix}   \begin{pmatrix}
1 & 0& 0  & 0 \\ 0 & g_1 &  0 & 0 \\ 0 &0 & g_2 & 0 \\  0 &0 & 0 & g_3
      \end{pmatrix}  k'\right) \right|\leq C \|(g_1,g_2,g_3)\|_{G_{2p}}^{-N}
      \end{align*}
      pour tout  $(g_1,g_2,g_3)\in [G_1\times G_p\times G_{p-1}  ]_{G_{2p}}$ et tous $n$ et $k'$ restant dans les compacts fixés plus haut. On en déduire la convergence et la majoration uniforme de l'intégrale sur $[G_1\times G_p\times G_{p-1}  ]_{G_{2p}}$ de la fonction ci-dessus multipliée par le facteur $|g_1|^{2p-1}  |\det(g_2)\det(g_3)|^{-1})^{s-1/2} $ pour tout $s\in \RR$. Finalement, il reste à vérifier la convergence de
      \begin{align*}
         \int_{\RR_+^\times}    \frac{t^{2ps+1}}{(  t^2+1)^{(p+1/2)(s+1/2)}} \, dt
      \end{align*}
 pour $s$ assez grand,     ce qui évident puisque pour $s>0$ on a 
 \begin{align*}
   &  \frac{t^{2ps+1}}{(  t^2+1)^{(p+1/2)(s+1/2)}} \sim_{t\to 0} t^{2ps+1}\\
       & \frac{t^{2ps+1}}{(  t^2+1)^{(p+1/2)(s+1/2)}} \sim_{t\to +\infty}   |t|^{1/2-s-p}.
 \end{align*}

 Une fois que la convergence absolue est vérifiée, on peut utiliser le calcul de $I_2(\la,\varphi)$ dans \eqref{eq:calcul-I2} qui fait apparaître l'intégrale
 \begin{align*}
   \int_{  [N_{1,2p-1}] } \varphi_\la\left( \begin{pmatrix} 1 & 0 \\ 0 & n  
      \end{pmatrix} g \right)\,dn
 \end{align*}
 qui est nulle puisque $\pi$ est cuspidale. Par conséquent, on a $I_2(\la,\varphi)=0$.
    \end{paragr}

  \subsection{Une contribution spectrale dans la formule des traces}\label{ssec:unecontrib}

  \begin{paragr}[Caractères relatifs.] --- \label{S:car rel} Soit $\bc_{\pi}$ une $K$-base de $\Ac_{P,\pi}(G)$ au sens de \cite[§ 2.8.3]{BCZ}. Rappelons qu'il s'agit d'une union, sur l'ensemble des classes d'isomorphisme des représentations irréductibles $\tau$ de $K$, de bases orthonormées, pour le produit de Petersson, du sous-espace de $\Ac_{P,\pi}(G)$ formé des fonctions qui se transforment sous $K$ selon $\tau$. Soit 
    \begin{align*}
      B:\Ac_{P,\pi}(G)\times \Ac_{P,\pi}(G)\to \CC
    \end{align*}
    une forme sesquilinéaire (à droite) et continue.

    On peut alors définir le caractère relatif: c'est la forme linéaire continue sur $\Sc(G(\AAA))$ qui dépend de $\la\in \ago_{P,\CC}^{G,*}$ et qui est donnée par
  \begin{align*}
    \Jc_B(\la,f)=\sum_{\varphi\in \bc_\pi} B(I_P(\la,f)\varphi,\varphi).
  \end{align*}
  La somme converge absolument,  uniformément pour $\Re(\la)$ dans un compact et définit une fonction holomorphe de la variable $\la$. De plus, la somme ne dépend pas du choix de  la $K$-base $\bc_{\pi}$. Pour toutes ces propriétés, on renvoie à  \cite[proposition  2.8.4.1]{BCZ}.
\end{paragr}

\begin{paragr} Selon la construction du § \ref{S:car rel}, on introduit le caractère relatif \og pondéré\fg,  associé à $(P,\pi)$:
\begin{align}\label{eq:car pondere}
    J_{P,\pi}(f)= \sum_{\varphi\in \bc_\pi}   J_P(\mc_{P,\pi} I_P(f)\varphi) \overline{J_P(\varphi)}.
\end{align}
C'est un analogue relatif du caractère pondéré introduit par Arthur, cf. \cite[section 15]{Ar-cours}.
\end{paragr}

\begin{paragr}[Résultat principal.] --- Nous sommes désormais en mesure d'énoncer le résultat principal de notre note, à savoir le calcul de la distribution $J_\chi$ en terme du caractère relatif pondéré qu'on vient d'introduire.
 
  \begin{theoreme}\label{thm:Jchi}
    Soit $\chi$ la donnée cuspidale associée à $(M,\pi)$. Pour tout $f\in \Sc(G(\AAA))$, on a
    \begin{align*}
      J_\chi(f)=   J_{P,\pi}(f).
    \end{align*}
  \end{theoreme}

  La démonstration de ce théorème sera donnée à la section \ref{ssec:demo-thmpal} et fera suite à une étude spectrale, à la section suivante, de l'intégrale sur $[H]$ de la $\chi$-composante du noyau de périodes de certaines séries d'Eisenstein. Avant de terminer cette section, nous aimerions donner des formules de covariance pour la distribution $J_{P,\pi}$  et montrer que celles-ci sont compatibles aux formules générales pour $J_\chi$ données dans \cite{CLi}.
\end{paragr}

\begin{paragr}[Covariance de $J_{P,\pi}$.] ---  \label{S:covJPpi}Pour tout $f\in \Sc(G(\AAA))$ et $g\in G(\AAA)$, on note  $^g\!f$ et $f^g$ les fonctions dans $\Sc(G(\AAA))$ définies par    $^g\!f(x)=f(gx)$ et $f^g(x)=f(xg)$ pour tout $x\in G(\AAA)$.

  \begin{proposition}
     \label{prop:cov-JPpi}
     Pour tout $y\in H(\AAA)$ et tout $f\in \Sc(G(\AAA))$, on a
    \begin{align*}
    J_{P,\pi}(f^y)&=J_{P,\pi}(f)\\
      J_{P,\pi}(\,  ^y\!f)- J_{P,\pi}(f)&=\sum_{\varphi\in \bc_\pi}   J_P(I_P'(y^{-1}) I_P(f)\varphi) \overline{J_P(\varphi)}\\
  &    +\sum_{\varphi\in \bc_{\pi'}}   J_Q(I_Q'(y^{-1})  I_Q(f)\varphi) \overline{J_Q(\varphi)},
    \end{align*}
    où $\bc_{\pi'}$ est une $K$-base de de $\Ac_{Q,\pi'}(G)$.
  \end{proposition}

  \begin{preuve}
    On a     $I_P(\,  ^y\!f)=  I_P(y^{-1}) I_P(f)$ et  $I_P(f^y)=  I_P(f)I_P(y^{-1}) $ pour tout $y\in G(\AAA)$. Prenons désormais $y\in H(\AAA)$ Par un changement de $K$-base, on a:
    \begin{align*}
      J_{P,\pi}(f^y)&= \sum_{\varphi\in \bc_\pi}   J_P(\mc_{P,\pi} I_P(f)I_P(y^{-1})\varphi) \overline{J_P(\varphi)}\\
      &=\sum_{\varphi\in \bc_\pi}   J_P(\mc_{P,\pi} I_P(f)\varphi) \overline{J_P(I_P(y)\varphi)}.
    \end{align*}
    Cette dernière expression n'est autre que $J_{P,\pi}(f)$ vu que la période d'entrelacement $J_P$ est $H(\AAA)$-invariante, cf. proposition \ref{prop:periode entrelac}. Cela donne la première égalité.

    Pour la seconde, on a, par le lemme \ref{lem:crochet},
    
    \begin{align*}
      J_{P,\pi}(\,  ^y\!f)= \sum_{\varphi\in \bc_\pi}   J_P(\mc_{P,\pi} I_P(y^{-1}) I_P(f)\varphi) \overline{J_P(\varphi)}\\= \sum_{\varphi\in \bc_\pi}   J_P(I_P(y^{-1}) \mc_{P,\pi}  I_P(f)\varphi) \overline{J_P(\varphi)}
      +\sum_{\varphi\in \bc_\pi}   J_P(I_P'(y^{-1}) I_P(f)\varphi) \overline{J_P(\varphi)}\\
        +\sum_{\varphi\in \bc_\pi}   J_P( M(w_0,0)^{-1} I_Q'(y^{-1}) M(w_0,0) I_P(f)\varphi) \overline{J_P(\varphi)}.
    \end{align*}
    En utilisant de nouveau l'invariance sous l'action de $H(\AAA)$ de la période d'entrelacement $J_P$, on voit qu'on a
    \begin{align*}
      \sum_{\varphi\in \bc_\pi}   J_P(I_P(y^{-1}) \mc_{P,\pi}  I_P(f)\varphi) \overline{J_P(\varphi)}= J_{P,\pi}(f).
    \end{align*}

    En utilisant l'équation fonctionnelle du corollaire \ref{cor:LaTEis} ainsi que le changement de base $\varphi\mapsto M(w_0,0)\varphi$ qui envoie $\bc_{\pi}$ sur une $K$-base $\bc_{\pi'}$ de $\Ac_{Q,\pi'}(G)$, on voit qu'on a
    \begin{align*}
      \sum_{\varphi\in \bc_\pi}   J_P( M(w_0,0)^{-1} I_Q'(y^{-1}) M(w_0,0) I_P(f)\varphi) \overline{J_P(\varphi)}\\
      =  \sum_{\varphi\in \bc_{\pi'}}   J_Q(I_Q'(y^{-1})  I_Q(f)\varphi) \overline{J_Q(\varphi)}.
           \end{align*}
    
  \end{preuve}
\end{paragr}

\begin{paragr}[Covariance de $J_\chi$.] --- Nous allons formuler dans notre situation particulière un résultat général \cite[proposition 3.5.5.1]{CLi}.  Soit $y\in H(\AAA)$ et $K^H=K_{p+1}\times K_p$ qui est un sous-groupe compact maximal de $H(\AAA)$.   On définit alors $f_{P,y}\in \Sc( M(\AAA))$ par la formule suivante: 
  \begin{align}\label{eq:fPy}
\forall m\in M(\AAA), \   f_{P,y} (m)=  \int_{K^H\times K^H}  \int_{N_P(\AAA)} f(  k_1^{-1}  m  n k_2) dn   \frac{\bg \al ,-H_P(k_1y)\bd}{\theta_P(\al)} dk_1dk_2.
  \end{align}

  Soit $K_{\chi,f_{P,y} }$ la composante, attachée à la donnée cuspidale dans $\Xgo(M)$ associée à $(M,\pi)$,  du noyau défini par:
  \begin{align*}
(x_1,x_2)\in [M]\times [M]\mapsto      \sum_{\ga\in M(F)}  f_{P,y}(x_1^{-1}   \ga x_2) .
  \end{align*}
  Notons que cette donnée cuspidale associée à $(M,\pi)$  est l'unique  antécédent de $\chi$  par l'application naturelle $\Xgo(M)\to \Xgo(G)$, ce qui justifie la notation.
  
  On vérifie que l'application
  \begin{align*}
  x_1\in  M^H(\AAA)\mapsto   \delta_{P^H}(x_1)^{-1} \int_{[M^H]}    \delta_{P^H}(x_2)^{-1} \exp(\bg 2\rho_P^G,H_P(x_2)K_{\chi,f_{P,y} }(x_1,x_2) \, dx_2
  \end{align*}
  est invariante par le centre $Z_M(\AAA)$ de $M(\AAA)$.
  On pose alors 
  \begin{align*}
    J_\chi^P( f_{P,y})=\int_{[M^H]_M}  \int_{[M^H]}  \delta_{P^H}(x_1)^{-1}  \delta_{P^H}(x_2)^{-1} \exp(\bg 2\rho_P^G,H_P(x_2)\bd) K_{\chi,f_{P,y} }(x_1,x_2) \, dx_2dx_1,
  \end{align*}
  où $\delta_{P^H}$ est le caractère modulaire du groupe $P^H(\AAA)$.     Notons que l'intégrale est ici absolument convergente car l'intégrande est à décroissance rapide.

  Rappelons qu'on a défini un élément  $w_3\in   \, _Q W_H$  au lemme \ref{lem:classes}.

Soit $f_{Q,y}\in \Sc( L(\AAA))$  définie par la formule:
\begin{align}\label{eq:fQy}
  \forall m\in L(\AAA), \   f_{Q,y} (m)=  \int_{K^H\times K^H}  \int_{N_Q(\AAA)} f(   (w_3k_1)^{-1} m  n w_3 k_2) dn   \frac{\bg \be,-H_Q(w_3k_1y)\bd}{\theta_Q(\al)} dk_1dk_2
  \end{align}
  où $w_3\in   \, _Q W_H$ est défini au lemme \ref{lem:classes}. Soit $K_{\chi,f_{Q,y} }$ la composante, attachée à la donnée cuspidale dans $\Xgo(L)$ définie par $(L,\pi')$,   du noyau défini par:
  \begin{align*}
(x_1,x_2)\in [L]\times [L]\mapsto     \sum_{\ga\in L(F)}  f_{Q,y}(x_1^{-1}   \ga x_2)\,dz.
  \end{align*}
Avec les notations du § \ref{S:JQ}, on pose:
  \begin{align*}
    J_\chi^Q( f_{Q,y})=\int_{[M^H_3]_{M_3}}  \int_{[M^H_3]}  \delta_{Q_3^H}(x_1)^{-1}  \delta_{Q_3^H}(x_2)^{-1} \exp(\bg 2\rho_{Q}^G,H_P(w_3x_2)\bd) K_{\chi,f_{Q,y} }(w_3x_1w_3^{-1},w_3x_2w_3^{-1}) \, dx_1dx_2.
  \end{align*}

  \begin{proposition}   \label{prop:covJchi}\cite[proposition 3.5.5.1]{CLi} Pour tout $y\in H(\AAA)$ et tout $f\in \Sc(G(\AAA))$, on a
    \begin{align*}
      J_\chi(f^y)&=J_\chi(f)\\
     J_\chi(  ^y\!f)- J_\chi(f)&=    J_\chi^P( f_{P,y})+ J_\chi^Q( f_{Q,y}).
    \end{align*}
  \end{proposition}

  \begin{preuve}
    La première égalité résulte immédiatement de l'assertion 1 de \cite[proposition 3.5.5.1]{CLi}. Pour la seconde égalité, on utilise l'assertion 2 de  \cite[proposition 3.5.5.1]{CLi}. Cette assertion écrit la différence $J_\chi(  ^y\!f)- J_\chi(f)$ comme une somme  de contributions  indexées par des triplets $(R,w_1',w_2')$ formés d'un sous-groupe parabolique standard $R\subsetneq G$ et d'éléments $w_i'\in \, _R W_H$ pour $i=1,2$ tels que
    \begin{align}\label{eq:flat}
          R\in \fc^{G,\flat}(P_0, w_1'\theta {w_1'}^{-1},w_2'\theta {w_2'}^{-1}),
    \end{align}
    avec les notation de \cite{CLi}.  Vu la forme particulière de la donnée cuspidale $\chi$ qu'on considère ici, seuls les sous-groupes paraboliques $R=P$ et  $R=Q$ peuvent donner des contributions non nulles.

    Considérons d'abord la contribution de $R=P$. La condition \eqref{eq:flat} est vérifiée si et seulement si $w_1',w_2'$ sont égaux  l'élément trivial de $W$. Dans ce cas, la fonction $f_{P,y}$ qu'on a définie en \eqref{eq:fPy}, est égale à la fonction définie en  \cite[éq. (3.5.4.1)]{CLi} pour les valeurs $Q=P, w_1'=w_2'=1$ et le caractère $\eta$ trivial. En effet, pour $w_1'=w_2'=1$, l'expression notée $p_{w_1',w_2'}^P(0,X)$ dans \cite{CLi} pour $X\in \ago_P^G$, est égale à
    \begin{align*}
      \int_{\ago_P^G} (\hat\tau_P(H)-\hat\tau_P(H-X))\, dH=   \frac{\bg \al,X\bd  }{ \theta_P(\al) }.
    \end{align*}
\emph{Stricto sensu}, il faudrait remplacer $\hat\tau_P(H)$ par  $\tau_P(H)$ avec $\tau_P$ la fonction caractéristique de la chambre obtuse   mais,  comme ici $P$ est maximal, on a $\tau_P = \hat\tau_P$.
 La contribution, associée dans \cite[proposition 3.5.5.1]{CLi} à $P$ et $w_1'=w_2'=1$, est le terme constant d'un polynôme-exponentielle en $T$ qui, vu la donnée cuspidale considérée ici, ne dépend pas de $T$. On en déduit  que la contribution associée est donnée par   $J_\chi^P( f_{P,y})$.

 Considérons ensuite la  contribution de $R=Q$. La condition \eqref{eq:flat} est vérifiée si et seulement si $w_1'=w_2'=w_3$. Cette fois-ci la fonction définie en  \cite[éq. (3.5.4.1)]{CLi}, pour   $w_1'=w_2'=w_3$ et le caractère $\eta$ trivial, est égale à  la fonction $f_{Q,y}$ définie en \eqref{eq:fQy}. En effet, l'expression $p_{w_3,w_3}^Q(0,X)$ se calcule comme ci-dessus et il suffit ensuite d'observer que, si l'on introduit $\theta_3=w_3\theta w_3^{-1}$ et  $H_3$ le centralisateur de $\theta_3$ , on a     $K\cap H_3=w_3 (K\cap H) w_3^{-1}$ et un changement de variables donne l'égalité cherchée. Un changement de variables et la même observation que ci-dessus (à savoir qu'un certain polynôme-exponentielle en $T$ est constant) identifie ensuite la contribution à  $J_\chi^Q( f_{Q,y})$.
\end{preuve}

Les formules données par les propositions \ref{prop:cov-JPpi} et \ref{prop:covJchi} sont, comme il se doit, compatibles comme le montre le lemme suivant.

\begin{lemme}
  \label{lem:compatibilite}
  On a
  \begin{align*}
    &    J_\chi^P( f_{P,y})=\sum_{\varphi\in \bc_\pi}   J_P(I_P'(y^{-1}) I_P(f)\varphi) \overline{J_P(\varphi)}\\
    &J_\chi^Q( f_{Q,y})=\sum_{\varphi\in \bc_{\pi'}}   J_Q(I_Q'(y^{-1})  I_Q(f)\varphi) \overline{J_Q(\varphi)}.
  \end{align*}
\end{lemme}

\begin{preuve}
  On donne la démonstration pour la première égalité, la seconde se traitant de manière analogue. Comme au § \ref{S:intro-DSFT}, on dispose du noyau $K_{P,\chi,f}$ attaché à $\chi$. En utilisant la définition de $f_{P,y}$  et un changement de variables, on voit qu'on a pour tous $x_1,x_2\in M^H(\AAA)$
  \begin{align*}
    K_{\chi,f_{P,y} }(x_1,x_2)=  \exp(-\bg 2\rho_P,H_P(x_2)\bd)  \int_{K^H\times K^H}  K_{P,\chi,f}(x_1k_1,x_2k_2)       \frac{\bg \al ,-H_P(k_1y)\bd}{\theta_P(\al)} dk_1dk_2.
  \end{align*}
  Il s'ensuit qu'on a, avec les notations des §§ \ref{S:JP}, \ref{S:M derive} et \ref{S:covJPpi}, en particulier avec \eqref{eq:fct-kP} 
  \begin{align*}
      J_\chi^P( f_{P,y})&=  \int_{ [H]_P\times [H]_P}     K_{P,\chi,f}^0(x_1,x_2)       \frac{\bg \al ,-H_P(k_P(x) y)\bd}{\theta_P(\al)} \, dx_1dx_2
  \end{align*}
  et
  \begin{align*}
    K_{P,\chi,f}^0(x_1,x_2)   =\int_{[Z_M]}  K_{P,\chi,f}(x_1,zx_2)\, dz.   
  \end{align*}
  On a $-H_P(k_P(x) y)= H_P(x) -H_P(xy)$ et $H_P(xy^{-1}) -H_P(x)=H_P(k_P(x) y^{-1})$.

  Par le changement de variables $x_1\mapsto xy_2$, il vient:
    \begin{align*}
      J_\chi^P( f_{P,y})&=  \int_{ [H]_P\times [H]_P}     K_{P,\chi,f}^0(x_1y^{-1},x_2)       \frac{\bg \al ,H_P(k_P(x_1) y^{-1})\bd}{\theta_P(\al)} \, dx_1dx_2.
    \end{align*}
On conclut facilement en utilisant l'écriture suivante du noyau:
    \begin{align*}
        K_{P,\chi,f}^0(x_1y^{-1},x_2)  = \sum_{\varphi \in \bc_{\pi}}      (I_P(f) \varphi)(x_1y^{-1})  \overline{\varphi(x_2)}.
    \end{align*}
\end{preuve}
\end{paragr}

\subsection{Une contribution spectrale pour l'espace symétrique}\label{ssec:undvptspec}
    
\begin{paragr} Voici le principal résultat de cette section.
  
  \begin{theoreme}\label{thm:chi dvpt}
    Pour toute fonction $f\in \Sc(G(\AAA))$ et tout $x\in  [G]_G$, on a 
    \begin{align}\label{eq:chi dvpt}
    \int_{[H]}K_{\chi,f}(x,y)\,dy = \sum_{\varphi\in \bc_\pi}   E(x,I_P(f)\varphi,0) \overline{J_P(\varphi) },
    \end{align}
    égalité dans laquelle l'intégrale et la somme sont absolument convergentes.

    En outre, l'intégrale ci-dessous est non nulle pour au moins une fonction $f$ et un élément $x\in  [G]_G$ si et seulement si la représentation $\sigma$ est symplectique.
  \end{theoreme}

  La dernière assertion du théorème résulte de \eqref{eq:chi dvpt} et de la proposition \ref{prop:periode entrelac} assertion 2. Notons que le terme constant  de  $E(x,I_P(f)\varphi,0)$ le long de $P$ est donné par \eqref{eq:cst-P} et qu'il est clairement non identiquement nul. Le reste de la sous-section est consacrée à la preuve  du développement \eqref{eq:chi dvpt}. Plus précisément, celui-ci résulte de la combinaison du lemme \ref{lem:lim T} et du lemme \ref{lem:calcul lim} ci-dessous.
\end{paragr}

\begin{paragr}[Une limite.] --- Désormais, on fixe $x\in G(\AAA)$. La convergence absolue du  membre de de \eqref{eq:chi dvpt} résulte de la continuité de $J_P$, cf. proposition \ref{prop:periode entrelac}, et celle de la série d'Eisenstein  ainsi que  des rappels du § \ref{S:car rel}. Par ailleurs, l'application $y\in [H]\mapsto K_{\chi,f}(x,y)$ est à décroissance rapide (cf. \cite[lemme 2.10.1.1]{BCZ}) d'où la convergence absolue de l'intégrale dans \eqref{eq:chi dvpt}.  Dans la suite, on fixe $f$ et on omet l'indice $f$ dans les notations.  On pose pour tout $y\in [H]_G$
    \begin{align*}
      K_\chi^0(x,y)= \int_{Z_G(\AAA)} K_\chi(x,ay)\, da.
    \end{align*}
Pour tout paramètre de troncature $T\in \ago_{P_0}^G$ assez positif,  on note alors $K_\chi^0\La^T_\theta$ le noyau $K_\chi^0$ lorsqu'on lui  applique l'opérateur de troncature relatif sur la deuxième variable.

    \begin{lemme} \label{lem:lim T}On a
    \begin{align*}
  \lim_{\bg \al,T\bd\to+\infty}          \int_{[H]_G} (K_\chi^0\La^T_\theta) (x,y)\,dy = \int_{[H]}K_\chi(x,y)\,dy .
    \end{align*}
    où l'intégrale de gauche converge absolument.
  \end{lemme}

  \begin{preuve}
    Notons que, à $x, y\in [H]_G$ fixés, on a, cf. remarque \ref{rq:op-zy},
    \begin{align*}
      \lim_{\bg \al,T\bd\to+\infty}        (K_\chi^0\La^T_\theta) (x,y)=  K_\chi^0(x,y).
    \end{align*}
    Le lemme est alors une conséquence immédiate du théorème de convergence dominée. Il reste donc à vérifier que ses hypothèses sont satisfaites. Rappelons qu'Arthur a introduit une fonction $g\mapsto F^G(\cdot,T)$ qui est la fonction caractéristique d'un compact de  $[G]_G$. On a alors la majoration:
    \begin{align*}
      |(K_\chi^0\La^T_\theta) (x,y)|&\leq |(K_\chi^0\La^T_\theta) (x,y)-  K_\chi^0(x,y)F^G(y,T)| + |K_\chi^0(x,y)F^G(y,T)|\\
      &\leq |(K_\chi^0\La^T_\theta) (x,y)-  K_\chi^0(x,y)F^G(y,T)| + |K_\chi^0(x,y)|.
    \end{align*}
    La fonction  $y\in [H]_G\mapsto |K_\chi^0(x,y)|$ est intégrable sur  $[H]_G$ car à décroissance rapide, cf. \cite[lemme 2.10.1.1]{BCZ}. Fixons $N>0$ de sorte que $y\mapsto \| y\|^{-N}_H$ est intégrable sur $[H]_G$. D'après \cite[proposition 2.4.3.2]{CLi}, il existe des éléments $X_1,\ldots,X_r$ dans l'algèbre enveloppante de l'algèbre de Lie complexifiée de $G$ tels que pour tout $y\in [H]_G$, on ait
    \begin{align*}
      |(K_{\chi,f}^0\La^T_\theta) (x,y)-  K_{\chi,f}^0(x,y)F^G(y,T)| \leq \|y\|_{H}^{-N} \sup_{1\leq i\leq r, g\in  [G]_G   } (\|g\|_G^{-N}| K^0_{\chi,R(X_i)f}(x,g)|).
    \end{align*}
    Ici $R$ désigne l'action régulière à droite de l'algèbre enveloppante sur l'espace de Schwartz. La borne supérieure est finie puisque $g\in  [G]_G \mapsto | K_{\chi,R(X_i)f}(x,g)|$ est à décroissance rapide, toujours par \cite[lemme 2.10.1.1]{BCZ}. On en déduit que les  hypothèses du théorème de convergence dominée sont satisfaites.
  \end{preuve}
\end{paragr}  

\begin{paragr}[Décomposition spectrale de la $\chi$-composante du noyau tronqué.] ---   La décomposition de Langlands de $K^0_\chi$ est donnée  pour tous $x,y\in [G]_G$ par:
  \begin{align}\label{eq:Langlands}
    K_\chi^0(x,y)=  \int_{i\ago_{P}^{G,*}} \sum_{\varphi\in \bc_\pi}   E(x,I_P(\la,f)\varphi,\la) \overline{ E(y,\varphi,\la) }\, d\la.
  \end{align}

  C'est une variante de la formule \cite[(12.7) p. 67]{Ar-cours}: a priori celle-ci comporte un autre terme associé à $(Q,\pi')$. En utilisant l'équation fonctionnelle des séries d'Eisenstein, on obtient bien la formule ci-dessus. Le développement spectral du noyau tronqué est donné par le lemme suivant.

  \begin{lemme}\label{lem:dvpt integ noy}
    Pour tous $x,y\in [H]_G$, on a 
  \begin{align*}
    (K_\chi^0\La^T_\theta) (x,y)=  \int_{i\ago_{P}^{G,*}} \sum_{\varphi\in \bc_\pi}   E(x,I_P(\la,f)\varphi,\la) \overline{ (\La^T_\theta E)(y,\varphi,\la) }\, d\la.
  \end{align*}
  \end{lemme}

  \begin{preuve} Par le théorème de Dixmier-Malliavin, on a  $\Sc(G(\AAA))=\Sc(G(\AAA))\ast\Sc(G(\AAA))$ où    $\ast$ est  le produit de convolution.   Sans perte de généralité, on peut donc  supposer qu'on a   $f=f_1\ast f^*_2$ avec  $f_1,f_2\in  \Sc(G(\AAA))$ et  $f_2^*(g)=\overline{f_2(g^{-1})}$ pour tout $g\in G(\AAA)$.  Un argument standard de changement de base montre alors que, pour tout $\la\in \ago_{P,\CC}^{G,*}$, on a
    
 \begin{align}
\label{eq:BC}   \sum_{\varphi\in \bc_\pi}   E(x,I_P(\la,f)\varphi,\la) \overline{ E(y,\varphi,\la) }=    \sum_{\varphi\in \bc_\pi}   E(x,I_P(\la,f_1)\varphi,\la) \overline{ E(y,I_P(-\bar \la,f_2)\varphi,-\bar \la) }\\
  \label{eq:BC2}     \sum_{\varphi\in \bc_\pi}   E(x,I_P(\la,f)\varphi,\la) \overline{ (\La^T_\theta E)(y,\varphi,\la) }= \sum_{\varphi\in \bc_\pi}   E(x,I_P(\la,f_1)\varphi,\la) \overline{ (\La^T_\theta E)(y,I_P(-\bar \la,f_2)\varphi,\la) }.
  \end{align}
Dans la suite, on prend $\la\in i\ago_{P}^{G,*}$ de sorte qu'on a $-\bar \la=\la$.
    
Soit $y\in H(\AAA)$. Vu la définition de $\La^T_\theta$ donnée en \eqref{eq:intro-LaTtheta} et le fait que les sommes dans  \eqref{eq:intro-LaTtheta} sont finies (à $x$ fixé),  les expressions $(K_\chi^0\La^T_\theta) (x,y)$ et $ (\La^T_\theta )E(y,I_P(\la,f_2)\varphi,\la)$ sont toutes deux données par une somme finie, indexée par le même ensemble fini de sous-groupes paraboliques standard $R$ de $G$ et d'éléments $\delta\in G(F)$,  de termes qui sont respectivement de la forme
    \begin{align*}
       \int_{[N_R]}   K_\chi^0(x,n \delta y)\,dn \ \ 
\text{        et   } \ \ 
      \int_{[N_R]}   E(n \delta y,I_P(\la,f_2)\varphi,\la)\,dn.
    \end{align*}
    On va prouver qu'on a
    \begin{align}
      \label{eq:pour Fubi}
    \int_{i\ago_{P}^{G,*}} \sum_{\varphi\in \bc_\pi}   \left| E(x,I_P(\la,f_1)\varphi,\la) \cdot \overline{ \int_{[N_R]} E(n\delta y,I_P(\la,f_2)\varphi,\la) \,dn }\right| \, d\la <\infty.
    \end{align}
    L'énoncé cherché résultera  alors d'une application du théorème de Fubini.    En utilisant l'inégalité de Cauchy-Schwartz, on doit que \eqref{eq:pour Fubi} résulte des deux majorations suivantes:
    \begin{align}
\label{eq:diag noy}      \int_{i\ago_{P}^{G,*}} \sum_{\varphi\in \bc_\pi}   \left| E(x,I_P(\la,f_1)\varphi,\la) \right|^2\, d\la =K^0_{\chi,f_1*f_1^*}(x,x) <\infty.
    \end{align}
    et
    \begin{align*}
      &          \int_{i\ago_{P}^{G,*}} \sum_{\varphi\in \bc_\pi}  \left|\int_{[N_R]} E(n\delta y,I_P(\la,f_2)\varphi,\la) \,dn\right|^2\, d\la \\
      &\leq  \int_{i\ago_{P}^{G,*}} \sum_{\varphi\in \bc_\pi}  \int_{[N_R]} \left| E(n\delta y,I_P(\la,f_2)\varphi,\la)   \right|^2 \,dn \, d\la\\
        &=  \int_{[N_R]} K_{\chi,f_2*f_2^*}^0( n\delta y,n\delta y) <\infty. 
    \end{align*}
 On notera que pour $i=1,2$, l'application $x\in G(\AAA)\mapsto K_{\chi,f_i*f_i^*}^0(x,x) $ est continue et positive.
\end{preuve}

\end{paragr}

\begin{paragr}[Décomposition spectrale de la période de la  $\chi$-composante.] ---   On introduit ensuite le  caractère relatif 
 \begin{align}
\label{eq:ecT}   \ec^T(x,f,\la) =\sum_{\varphi\in \bc_\pi}   E(x,I_P(\la,f)\varphi,\la) \overline{ \int_{[H]_G} (\La^T_\theta E)(y,\varphi,-\bar \la) \,dy}
 \end{align}
pour $\la\in \ago_{P,\CC}^{G,*}$.

  \begin{lemme}\label{lem:dvpt integ noy2}
    Pour  tout $x\in  [G]_G$, on a
     \begin{align*}
    \int_{[H]_G}(K_{\chi,f}^0\La^T_\theta)(x,y)\,dy = \int_{i\ago_{P}^{G,*}}  \ec^T(x,f,\la)  \, d\la.
    \end{align*}
  \end{lemme}

  \begin{preuve}  Comme dans la preuve du lemme \ref{lem:dvpt integ noy}, on peut et on va supposer qu'on a $f=f_1\ast f^*_2$ avec  $f_1,f_2\in  \Sc(G(\AAA))$. On a alors:
    \begin{align*}
      \ec^T(x,f,\la)  = \sum_{\varphi\in \bc_\pi}   E(x,I_P(\la,f_1)\varphi,\la) \overline{  \int_{[H]_G}(\La^T_\theta E)(y,I_P(-\bar \la,f_2)\varphi,\la) \, dy}.
    \end{align*}
    Le lemme résulte alors immédiatement du lemme \ref{lem:dvpt integ noy} couplé avec l'égalité  \eqref{eq:BC} et du théorème de Fubini puisque ce dernier permet d'échanger l'ordre d'intégration et de sommation.     Il s'agit donc de vérifier que les hypothèses du théorème de Fubini sont ici satisfaites, à savoir qu'on a
    \begin{align*}
   \int_{i\ago_{P}^{G,*}}   \sum_{\varphi\in \bc_\pi}   \left( \left|E(x,I_P(\la,f_1)\varphi,\la) \right| \cdot \int_{[H]_G} \left|(\La^T_\theta E)(y,I_P(\la,f_2)\varphi,\la)\right| \,dy \right)\, d\la<\infty.
    \end{align*}
    Par l'inégalité de Cauchy-Schwartz, il suffit prouver la finitude des deux expressions suivantes:
    \begin{align*}
      \int_{i\ago_{P}^{G,*}} \sum_{\varphi\in \bc_\pi}   \left| E(x,I_P(\la,f_1)\varphi,\la) \right|^2\, d\la 
    \end{align*}
    et
     \begin{align*}
      \int_{i\ago_{P}^{G,*}} \sum_{\varphi\in \bc_\pi}   \left( \int_{[H]_G} \left|(\La^T_\theta E)(y,I_P(\la,f_2)\varphi,\la)\right| \,dy\right)^2\, d\la .
     \end{align*}
     Pour la première expression, cela a été observé en \eqref{eq:diag noy}. Pour la seconde, on fixe $N>0$ assez grand pour que
     \begin{align*}
       v_N=\int_{[H]_G} \|g\|_G^{-2N}\,dg<\infty.
     \end{align*}
     Une nouvelle application de Cauchy-Schwartz, cette fois-ci à l'intégrale sur  $[H]_G$, nous ramène à majorer l'expression suivante:
     \begin{align*}
       v_N      \int_{i\ago_{P}^{G,*}} \sum_{\varphi\in \bc_\pi}  \int_{[H]_G}\|y\|_G^{2N} \left|(\La^T_\theta E)(y,I_P(\la,f_2)\varphi,\la)\right|^2 \,dy\, d\la\\
       =v_N  \int_{[H]_G}\|y\|_G^{2N}  \left(\int_{i\ago_{P}^{G,*}} \sum_{\varphi\in \bc_\pi} \left|(\La^T_\theta E)(y,I_P(\la,f_2)\varphi,\la)\right|^2 \, d\la\right)\,dy.
     \end{align*}
     En suivant la méthode de la preuve du lemme \ref{lem:dvpt integ noy}, on montre que l'expression entre parenthèses n'est autre que $(\La^T_\theta K_{\chi,f_2*f_2^*}^0\La^T_\theta)(y,y) $. Expliquons la notation. On obtient $\La^T_\theta K_{\chi,f_2*f_2^*}^0\La^T_\theta$ en appliquant l'opérateur de troncature  $\La^T_\theta $ à la variable de gauche de $K_{\chi,f_2*f_2^*}^0\La^T_\theta$. Pour conclure, il suffit d'observer que, pour tout $N>0$, il existe $C>0$ telle que pour tout $y\in [H]_G$ on ait
     \begin{align*}
         \left|(\La^T_\theta K_{\chi,f_2*f_2^*}^0\La^T_\theta )(y,y)\right| \leq C \|y\|^{-N}_G.
     \end{align*}

     Cela résulte de majorations du noyau $K_{\chi,f}$ de \cite[lemme 2.10.1.1]{BCZ} et des propriétés de l'opérateur  $\La^T_{\theta}$, cf. \cite[proposition 2.4.3.1]{CLi}.
   \end{preuve}
 \end{paragr}

 \begin{paragr}[Analyse asymptotique de caractères relatifs.] --- En complément du caractère relatif $\ec^T(x,f,\la)$ introduit en  \eqref{eq:ecT}, on introduit les caractères relatifs, pour $\la\in \ago_{P,\CC}^{G,*}$,
   \begin{align}
\label{eq:ecP}         &    \ec_P(x,f,\la) =\sum_{\varphi\in \bc_\pi}   E(x,I_P(\la,f)\varphi,\la) \overline{J_P(\varphi) };\\
  \label{eq:ecPb}    &\ec_{\bar P}(x,f,\la) =\sum_{\varphi\in \bc_\pi}   E(x,I_P(\la,f)\varphi,\la) \overline{J_{\bar P}(-\bar \la, \varphi) }.
  \end{align}

  \begin{lemme}\label{lem:car Schw}
    Les caractères relatifs $\ec^T(x,f,\la),   \ec_P(x,f,\la)$ et $ \ec_{\bar P}(x,f,\la) $ appartiennent à la classe de Schwartz de $  i\ago_{P}^{G,*}$  au sens où ils définissent des fonctions de $\la\in  i\ago_{P}^{G,*}$ qui sont  lisses et dont toutes les dérivées sont à décroissance rapide.
  \end{lemme}

  \begin{preuve}     Les trois caractères relatifs sont holomorphes  au voisinage de  $i\ago_{P}^{G,*}$ donc lisses sur cette droite. Notons que, pour tout $\la\in   i\ago_{P}^{G,*}$, ils seront reliés par la formule suivante, conséquence  du théorème \ref{thm:LaTEis}
    \begin{align}\label{eq:lien}
     \ec^T(x,f,\la)   = \frac{   -\ec_P(x,f,\la) \exp(-\bg \la,T_P\bd)+ \ec_{\bar P}(x,f,\la) \exp(-\bg \la,T_{\bar P}\bd)}{  \theta_P(\la) }.
          \end{align}
    
       Il suffit donc de prouver que   $ \ec_P(x,f,\la)$ et $ \ec_{\bar P}(x,f,\la) $ sont dans la classe de Schwartz. On va considérer  uniquement  $ \ec_{\bar P}(x,f,\la) $; le cas de $ \ec_P(x,f,\la)$ se traite de la même façon et est, en fait, sensiblement plus simple. Pour cela, on va reprendre  une méthode due à Lapid, cf. \cite[sections 5 et 6]{LapFRTF} et \cite[sections 3 et 4]{LapHC}. On va d'abord prouver la décroissance  de $\la\mapsto \ec_{\bar P}(x,f,\la)$ sur certains voisinages ouverts de  $i \ago_{P}^{G,*}$ de la forme
           \begin{align*}
    \rc_{c,l}=\{\la\in \ago_{P,\CC}^{G,*}\mid |\Re(\bg \la,\al^\vee\bd)| < c(1+|\Im(\bg \la,\al^\vee\bd)|)^{-l}\}
  \end{align*}
  dépendant de réels $c,l>0$. La décroissance de ses dérivées sur $i \ago_{P}^{G,*}$ en résultera immédiatement par la formule de Cauchy. Comme dans la preuve du lemme \ref{lem:dvpt integ noy}, on peut et on va supposer, sans perte de généralité, qu'on a  $f=f_1\ast f^*_2$ avec  $f_1,f_2\in  \Sc(G(\AAA))$. Alors on a:
  
  \begin{align}\label{eq:variante JP}
    \ec_{\bar P}(x,f,\la) =\sum_{\varphi\in \bc_\pi}   E(x,I_P(\la,f_1)\varphi,\la) \overline{J_{\bar P}( -\bar \la, I_P(-\bar\la,f_2)\varphi) }.
  \end{align}

    On écrit $K=K_\infty K^\infty$ avec $K_\infty=\prod_{v\in V_{F,\infty}} K_v$ et  $K^\infty=\prod_{v\in V_{F}^{\infty}} K_v$. On fixe un sous-groupe ouvert, compact et normal $J$ de $K^\infty$ tel que $f_1$ et $f_2$ soient bi-invariantes par $J$. Soit $\hat K _\infty$ l'ensemble des classes d'isomorphisme  des  représentations irréductibles unitaires de  $K_\infty$.  Pour tout $\tau\in \hat K_\infty$, soit $\Ac_{P,\pi}(G)^{\tau,J}\subset \Ac_{P,\pi}(G)$ le sous-espace de dimension finie des fonctions $J$-invariantes à droites et qui se transforment selon $\tau$ par translation à droite sous $K_\infty$. Essentiellement par définition d'une $K$-base $\bc_\pi$,  un élément $\varphi\in  \bc_\pi$ contribue non trivialement seulement s'il appartient à un sous-espace $\Ac_{P,\pi}(G)^{\tau,J}$ pour un certain $\tau\in \hat K_\infty$. Pour tout $\tau\in \hat K_\infty$, on fixe alors une base $\bc_\pi(\tau,J)$ de $\Ac_{P,\pi}(G)^{\tau,J}$ pour le produit de Petersson. On peut et on va remplacer dans \eqref{eq:variante JP}  la $K$-base $\bc_\pi$ par la famille $\bc_\pi(J)=\cup_{\tau\in \hat K_\infty} \bc_\pi(\tau,J)$. On associe à chaque représentation $\tau\in \hat K_\infty$ une mesure $e_\tau$ supportée sur $K_\infty$, qui est un idempotent pour le produit de convolution et qui vérifie, pour tout  $\varphi\in \Ac_{P,\pi}(G)$   qui est $J$-invariant à droite,  $I_{P}(\la,e_{\tau})\varphi=\varphi$ si et seulement si  $\varphi\in \Ac_{P,\pi}(G)^{\tau,J}$. %
    On a alors
    \begin{align*}
      \ec_{\bar P}(x,f,\la) =\sum_{\tau,\tau_1,\tau_2\in  \hat K_\infty}\sum_{\varphi\in \bc_\pi(\tau,J)}   E(x,I_P(\la,e_{\tau_1}*f_1)\varphi,\la) \overline{J_{\bar P}( -\bar \la, I_P(-\bar\la,   e_{\tau_2}*f_2)\varphi) }.
    \end{align*}
    
    En utilisant l'inégalité triangulaire puis l'inégalité de Cauchy-Schwartz pour la somme sur $\tau$ et $\varphi$, on majore   $|\ec_{\bar P}(x,f,\la) | $ par le produit de

    \begin{align}
\label{eq:dec rap}&      \sum_{\tau ,\tau_1\in  \hat K_\infty}\sum_{\varphi\in \bc_\pi(\tau,J)}   \left|E(x,I_P(\la,e_{\tau_1}* f_1)\varphi,\la)\right|^2 \\
 \label{eq:dec rap2} &    \sum_{\tau ,\tau_1\in  \hat K_\infty} \sum_{\varphi\in \bc_\pi(\tau,J)} \left|  J_{\bar P}(\la, I_P(\la,e_{\tau_1}* f_2)\varphi)  \right|^2.
    \end{align}

Il résulte alors de \cite[proposition 6.1]{LapFRTF} (on suit en fait les notations de \cite[théorème 3.9.2.1]{chaudouardsymmetric} qui est une généralisation de l'énoncé de Lapid au cas de représentations discrètes) qu'il existe $l>0$ tel que pour tout $q>0 $ il existe $c>0$ et une semi-norme continue $\| \cdot\|$  sur $\Sc(G(\AAA))$ tels que pour tout $\la\in \rc_{c,l}$ et toute fonction $f_1\in \Sc(G(\AAA))$  bi-invariante par $J$, on ait:
\begin{align*}
 \sum_{\tau,\tau_1 \in  \hat K_\infty}\sum_{\varphi\in \bc_\pi(\tau,J)}    |E(x,I_P(\la,e_{\tau_1}*f_1)\varphi,\la)|^2 \leq \frac{\| f_1\|^2}{(1+|\bg \la, \al^\vee\bd|^2)^q}.
\end{align*}

  L'ingrédient principal pour majorer \eqref{eq:dec rap2}  est le contrôle de la norme de l'opérateur d'entrelacement $M(w_0,\la)$. Pour tout $\tau\in \hat K_\infty$, on note $\la_\tau\geq 0$ la valeur propre de l'opérateur de Casimir. D'après \cite[corollaire 3.13]{LapHC} (là encore, on suit en fait les notations de \cite[proposition 3.2.4.1]{chaudouardsymmetric} qui généralise de l'énoncé de Lapid au cas de représentations discrètes), il existe $k,c,l,C$ des réels $>0$ tels que pour tout  $\tau\in \hat K_\infty$, tout $\varphi\in \Ac_{P,\pi}(G)^{\tau,J}$ et $\la\in  \rc_{c,l}$
  \begin{align*}
    \|M(w_0,\la)\varphi\|_Q\leq  C ( 1+ |\bg \la, \al^\vee\bd|^2+\la_\tau^2)^k  \|\varphi\|_P.
  \end{align*}

   Rappelons qu'on a $J_{\bar P}(\la,\varphi)=J_Q(M(w_0,\la)\varphi)$. En utilisant la continuité de $J_Q$ et  la majoration ci-dessus  (les  arguments sont ceux de la preuve de \cite[proposition 3.2.4.1]{chaudouardsymmetric}), on voit qu'il existe $l>0$ tel que, pour tout $q>0$, il existe $c>0$ et une semi-norme continue $\| \cdot\|$  sur $\Sc(G(\AAA))$ tels que pour tous $\tau,\tau_1\in\hat K_\infty$ et tout $\la\in \rc_{c,l}$  
  \begin{align*}
    \sum_{\varphi\in \bc_\pi(\tau,J)} \left|  J_{\bar P}(\la, I_P(\la,e_{\tau_1}*f_2)\varphi)  \right|^2 \leq \frac{ \| f_2 \|^2   }{(1+|\bg \la, \al^\vee\bd|^2)^q(1+\la_{\tau_1}^2)^{q} (1+\la_{\tau}^2)^{q}   }.
  \end{align*}
  Cela conclut la démonstration vu que, pour $q$ assez grand, on a 
  \begin{align*}
    \sum_{\tau \in\hat K_\infty} (1+\la_{\tau}^2)^{q}   <\infty.
  \end{align*}

 \end{preuve}
\end{paragr}

\begin{paragr}[Calcul d'une limite.] ---

  \begin{lemme}\label{lem:calcul lim}
    Pour tout $x\in [G]_G$, on a 
    \begin{align*}
      \lim_{\bg \al,T\bd\to+\infty}          \int_{[H]_G} (K_\chi^0\La^T_\theta) (x,y)\,dy = \ec_{P}(x,f,0).
    \end{align*}
    où le second membre est défini en \eqref{eq:ecP}.
  \end{lemme}
  
  \begin{preuve}
    En utilisant le lemme \ref{lem:dvpt integ noy2} ainsi que l'égalité \eqref{eq:lien}, il s'agit de calculer la limite pour $\bg \al,T\bd\to+\infty$ de la somme de
    \begin{align}\label{eq:integ 1}
      \int_{i\ago_P^{G,*}} \ec_{\bar P}(x,f,\la) \frac{  \exp(-\bg \la,T_{\bar P}\bd) -\exp(-\bg \la,T_P\bd)  }{  \theta_P(\la) }\, d\la
    \end{align}
    et
     \begin{align*}
      \int_{i\ago_P^{G,*}}\exp(-\bg \la,T_P\bd) \frac{ \ec_{\bar P}(x,f,\la)   -\ec_P(x,f,\la) }{  \theta_P(\la)} \, d\la.
     \end{align*}
     Cette dernière tend vers $0$ par le lemme de Riemann-Lebesgue vu que $\la\mapsto \frac{ \ec_{\bar P}(x,f,\la)   -\ec_P(x,f,\la) }{  \theta_P(\la)} $ est lisse à décroissance rapide, cf. lemme \ref{lem:car Schw}. Pour traiter l'intégrale \eqref{eq:integ 1}, on remarque qu'on a 
     \begin{align*}
       \frac{  \exp(-\bg \la,T_{\bar P}\bd) -\exp(-\bg \la,T_P\bd)  }{  \theta_P(\la) }= \int_{\ago_P^{G}} \phi^T(H) \exp(-\bg \la,H \bd)\, dH
     \end{align*}
c'est-à-dire que le membre de gauche est la transformée de Fourier de la fonction caractéristique  $\phi^T$ du segment $[T_{\bar P}; T_P]$. À l'aide de la formule de Plancherel, l'intégrale \eqref{eq:integ 1} est égale à
     \begin{align*}
       \int_{ \ago_P^{G} }  \phi^T(H) \psi(H)\,dH\to_{  \bg \al,T\bd\to+\infty}  \int_{ \ago_P^{G} }   \psi(H)\,dH
     \end{align*}
     où $\psi$ est la transformée de Fourier inverse de la fonction  $\la\mapsto \ec_{\bar P}(x,f,\la) $. En particulier l'intégrale de $\psi$ est la valeur en $\la=0$ de $ \ec_{\bar P}(x,f,\la) $. L'équation fonctionnelle des périodes d'entrelacement, cf. corollaire \ref{cor:LaTEis}, entraîne qu'on a
     \begin{align*}
       \ec_{\bar P}(x,f,0) =\ec_{P}(x,f,0). 
     \end{align*}
     Cela conclut la démonstration.
  \end{preuve}
\end{paragr}

\subsection{Démonstration du théorème \ref{thm:Jchi}}\label{ssec:demo-thmpal}

\begin{paragr} On a 
    \begin{align*}
    \int_{[H]}\La^T_\theta K_\chi(x,y)\,dy = \sum_{\varphi\in \bc_\pi}   (\La^T_\theta E)(x,I_P(\la,f)\varphi,0) \overline{J_P(\varphi) }.
    \end{align*}

    \begin{lemme}
      L'expression
      \begin{align}\label{eq:double integ}
         \int_{[H]_G}  \int_{[H]}\La^T_\theta K_\chi(x,y)\,dy dx 
      \end{align}
      comme fonction de $T$ est une fonction affine dont la partie constante est le caractère relatif pondéré $J_{P,\pi}(f)$ défini en \eqref{eq:car pondere}.
    \end{lemme}

    \begin{preuve}
          Il résulte du théorème \ref{thm:chi dvpt} que l'expression \eqref{eq:double integ} est égale à  
    \begin{align*}
    \sum_{\varphi\in \bc_\pi}      \left(\int_{[H]_G} \La^T_\theta E(x,I_P(f) \varphi,0)  dx\right)   \overline{J_P(\varphi) }.
    \end{align*}
    Le lemme est alors une conséquence immédiate du calcul de la série d'Eisenstein tronquée effectué au corollaire \ref{cor:LaTEisven 0}.
  \end{preuve}
\end{paragr}

  \bibliography{bibliographie}

\begin{thebibliography}{BPCZ22}

\bibitem[Art80]{ar2}
J.~Arthur.
\newblock A trace formula for reductive groups. {II}. {A}pplications of a
  truncation operator.
\newblock {\em Compositio Math.}, 40(1):87--121, 1980.

\bibitem[Art81]{arthur2}
J.~Arthur.
\newblock The trace formula in invariant form.
\newblock {\em Ann. of Math. (2)}, 114(1):1--74, 1981.

\bibitem[Art05]{Ar-cours}
J.~Arthur.
\newblock An introduction to the trace formula.
\newblock In {\em Harmonic analysis, the trace formula, and {S}himura
  varieties}, volume~4 of {\em Clay Math. Proc.}, pages 1--263. Amer. Math.
  Soc., Providence, RI, 2005.

\bibitem[Bad08]{Badu}
A.~I. Badulescu.
\newblock Global {J}acquet-{L}anglands correspondence, multiplicity one and
  classification of automorphic representations.
\newblock {\em Invent. Math.}, 172(2):383--438, 2008.
\newblock With an appendix by Neven Grbac.

\bibitem[{Beu}21]{RBP}
R.~{Beuzart-Plessis}.
\newblock {Comparison of local spherical characters and the Ichino-Ikeda
  conjecture for unitary groups}.
\newblock {\em J. Inst. Math. Jussieu}, 20(6):1803--1854, 2021.

\bibitem[BPCZ22]{BCZ}
R.~Beuzart-Plessis, P.-H. Chaudouard, and M.~Zydor.
\newblock The global {G}an-{G}ross-{P}rasad conjecture for unitary groups: the
  endoscopic case.
\newblock {\em Publ. Math. Inst. Hautes \'{E}tudes Sci.}, 135:183--336, 2022.

\bibitem[BR10]{BaduR}
A.~I. Badulescu and D.~Renard.
\newblock Unitary dual of {${\rm GL}(n)$} at {A}rchimedean places and global
  {J}acquet-{L}anglands correspondence.
\newblock {\em Compos. Math.}, 146(5):1115--1164, 2010.

\bibitem[Cha25]{chaudouardsymmetric}
P.-H. Chaudouard.
\newblock A spectral expansion for the symmetric space {${\rm GL}_n(E)/{\rm
  GL}_n(F)$}.
\newblock {\em Selecta Math. (N.S.)}, 31(4):Paper No. 71, 2025.

\bibitem[CL]{CLi}
P.-H. Chaudouard and H.~{Li}.
\newblock Sur la formule des traces de {G}uo-{J}acquet.

\bibitem[FJ93]{FriJa}
S.~Friedberg and H.~Jacquet.
\newblock Linear periods.
\newblock {\em J. Reine Angew. Math.}, 443:91--139, 1993.

\bibitem[Jac97]{Jac-Edin}
H.~Jacquet.
\newblock Automorphic spectrum of symmetric spaces.
\newblock In {\em Representation theory and automorphic forms ({E}dinburgh,
  1996)}, volume~61 of {\em Proc. Sympos. Pure Math.}, pages 443--455. Amer.
  Math. Soc., Providence, RI, 1997.

\bibitem[JLR99]{JLR}
H.~Jacquet, E.~Lapid, and J.~Rogawski.
\newblock Periods of automorphic forms.
\newblock {\em J. Amer. Math. Soc.}, 12(1):173--240, 1999.

\bibitem[Lap06]{LapFRTF}
E.~Lapid.
\newblock On the fine spectral expansion of {J}acquet's relative trace formula.
\newblock {\em J. Inst. Math. Jussieu}, 5(2):263--308, 2006.

\bibitem[Lap08]{Lap-remark}
E.~Lapid.
\newblock A remark on {E}isenstein series.
\newblock In {\em Eisenstein series and applications}, volume 258 of {\em
  Progr. Math.}, pages 239--249. Birkh\"{a}user Boston, Boston, MA, 2008.

\bibitem[Lap13]{LapHC}
E.~Lapid.
\newblock On the {H}arish-{C}handra {S}chwartz space of {$G(F)\backslash G(\Bbb
  A)$}.
\newblock In {\em Automorphic representations and {$L$}-functions}, volume~22
  of {\em Tata Inst. Fundam. Res. Stud. Math.}, pages 335--377. Tata Inst.
  Fund. Res., Mumbai, 2013.
\newblock With an appendix by Farrell Brumley.

\bibitem[MOY25]{matringe2025}
N.~Matringe, O.~Offen, and C.~Yang.
\newblock Intertwining periods, {L}-functions and local-global principles for
  distinction of automorphic representations, 2025.

\bibitem[MW94]{MWlivre}
C.~M{\oe}glin and J.-L. Waldspurger.
\newblock {\em D\'{e}composition spectrale et s\'{e}ries d'{E}isenstein},
  volume 113 of {\em Progress in Mathematics}.
\newblock Birkh\"{a}user Verlag, Basel, 1994.
\newblock Une paraphrase de l'\'{E}criture.

\bibitem[Sak13]{Sak-spherical}
Y.~Sakellaridis.
\newblock Spherical functions on spherical varieties.
\newblock {\em Amer. J. Math.}, 135(5):1291--1381, 2013.

\bibitem[SV17]{SaVen}
Y.~Sakellaridis and A.~Venkatesh.
\newblock Periods and harmonic analysis on spherical varieties.
\newblock {\em Ast\'{e}risque}, (396):viii+360, 2017.

\bibitem[Tak23]{takeda2023}
S.~Takeda.
\newblock On dual groups of symmetric varieties and distinguished
  representations of $p$-adic groups, 2023.

\bibitem[Zha15]{Zha2}
C.~Zhang.
\newblock On linear periods.
\newblock {\em Math. Z.}, 279(1-2):61--84, 2015.

\bibitem[Zyd22]{Zydor}
M.~Zydor.
\newblock Periods of automorphic forms over reductive subgroups.
\newblock {\em Ann. Sci. \'{E}c. Norm. Sup\'{e}r. (4)}, 55(1):141--183, 2022.

\end{thebibliography}
\bibliographystyle{alpha}

\begin{flushleft}
Pierre-Henri Chaudouard \\
Université Paris Cité\\
CNRS \\
IMJ-PRG \\
Bâtiment Sophie Germain\\
8 place Aurélie Nemours\\
F-75013 PARIS \\
France\\
\medskip
Institut Universitaire de France\\
\medskip

email:\\
Pierre-Henri.Chaudouard@imj-prg.fr \\
\end{flushleft}

\end{document}